\documentclass[11pt]{article}
\usepackage{amsfonts}
\usepackage{amssymb}
\usepackage{graphicx}
\usepackage{amsmath}
\usepackage{amsthm}
\usepackage{enumitem}
\usepackage{tikz}
\usepackage{cancel}
\usepackage{todonotes}
\usetikzlibrary{matrix}
\usetikzlibrary{arrows,shapes}
\usepackage{cite}
\usepackage{hyperref}
\hypersetup{colorlinks=true, urlcolor=blue}
\expandafter\let\expandafter\oldproof\csname\string\proof\endcsname
\let\oldendproof\endproof
\renewenvironment{proof}[1][\proofname]{
\oldproof[\ttfamily\scshape\bf #1.]
}{\oldendproof}
\def\O{{\cal O}}

\def\B{\mathbb{B}}
\def\R{{\rm I\!R}}
\def\N{{\rm I\! N}}
\def\ox{\bar{x}}
\def\oy{\bar{y}}

\def\ov{\bar{v}}
\def\ow{\bar{w}}

\def\op{\bar{p}}

\def\ve{\varepsilon}

\def\tilde{\widetilde}
\def\emp{\emptyset}

\def\Lm{{\Lambda}}
\def\tto{\rightrightarrows}

\def\d{{\rm d}}
\def\sub{\partial}

\def\disp{\displaystyle}

\def\Bar{\overline}
\def\ra{\rangle}
\def\la{\langle}
\def\ve{\varepsilon}
\def\N{{\rm I\!\!N}}

\def\ri{\mbox{\rm ri}\,}

\def\int{\mbox{\rm int}\,}
\def\gph{\mbox{\rm gph}\,}
\def\epi{\mbox{\rm epi}\,}

\def\dom{\mbox{\rm dom}\,}

\def\ker{\mbox{\rm ker}\,}

\def\cl{\mbox{\rm cl}\,}

\def\dn{\downarrow}
\def\O{\Omega}

\def\ph{\varphi}
\def\tph{\tilde\varphi}
\def\tPh{\tilde f}
\def\emp{\emptyset}
\def\st{\stackrel}
\def\oR{\Bar{\R}}
\def\lm{\lambda}
\def\olm{\bar\lambda}

\def\dd{\delta}
\def\al{\alpha}
\def\Th{\Theta}

\def\sm{\hbox{${1\over 2}$}}
\def\h{\hfill\Box}

\def\sce{\setcounter{equation}{0}}

\begin{document}
\vspace*{0.5in}
\begin{center}
{\bf STABILITY OF KKT SYSTEMS AND SUPERLINEAR CONVERGENCE OF THE SQP METHOD UNDER PARABOLIC REGULARITY}\\[2ex]
ASHKAN MOHAMMADI\footnote{Department of Mathematics, Wayne State University, Detroit, MI 48202, USA (ashkan.mohammadi@wayne.edu). Research of this author was partly supported by the US National Science Foundation under grant DMS-1808978 and by the US Air Force Office of Scientific Research under grant \#15RT0462.}
BORIS S. MORDUKHOVICH\footnote{Department of Mathematics, Wayne State University, Detroit, MI 48202, USA (boris@math.wayne.edu). Research of this author was partly supported by the US National Science Foundation under grants DMS-1512846 and DMS-1808978, by the US Air Force Office of Scientific Research under grant \#15RT0462, and by the Australian Research Council Discovery Project DP-190100555.} and M. EBRAHIM SARABI\footnote{Department of Mathematics, Miami University, Oxford, OH 45065, USA (sarabim@miamioh.edu).}
\end{center}
\vspace*{0.05in}
\small{\bf Abstract.} This paper pursues a two-fold goal. Firstly, we aim to derive novel second-order characterizations of important robust stability properties of perturbed Karush-Kuhn-Tucker systems for a broad class of constrained optimization problems generated by parabolically regular sets. Secondly, the obtained characterizations are applied to establish well-posedness and superlinear convergence of the basic sequential quadratic programming method to solve parabolically regular constrained optimization problems.\\[1ex]

\newtheorem{Theorem}{Theorem}[section]
\newtheorem{Proposition}[Theorem]{Proposition}
\newtheorem{Remark}[Theorem]{Remark}
\newtheorem{Lemma}[Theorem]{Lemma}
\newtheorem{Corollary}[Theorem]{Corollary}
\newtheorem{Definition}[Theorem]{Definition}
\newtheorem{Example}[Theorem]{Example}
\newtheorem{Algorithm}[Theorem]{Algorithm}
\renewcommand{\theequation}{{\thesection}.\arabic{equation}}
\renewcommand{\thefootnote}{\fnsymbol{footnote}}

\normalsize
\section{Introduction}\label{intro}

This paper is devoted to the second-order variational analysis and numerical applications for a general class of constrained optimization problems formulated in Section~\ref{sect2}. It has been well recognized in optimization theory and its applications that second-order analysis concerning both qualitative/theoretical and quantitative/numerical aspects of constrained optimization requires certain second-order regularity conditions. In this paper we explore a novel one, which is the {\em parabolic regularity} of the constraint set at the point in question.

Actually the notion of parabolic regularity for extended-real-valued functions was first formulated by Rockafellar and Wets \cite[Definition~13.65]{rw}, but was not studied there except for fully amenable compositions with the immediate application to the unconstrained format of second-order optimality; see also \cite[Subsection~3.3.5]{bs}. We are not familiar with any other publications where parabolic regularity was either further investigated, or applied to structural problems of constrained optimization.

Very recently \cite{mms19}, a systematic study of parabolically regular sets (both convex and nonconvex) has been conducted in our paper, where we reveal a fundamental role of this concept in second-order variational analysis and its important applications to optimization theory. In particular, it is shown in \cite{mms19} that parabolic regularity is {\em more general} than known second-order regularity notions, is {\em preserved} under various operations on sets, ensures---among other significant results---{\em twice epi-differentiability} of set indicator functions, {\em precise/equality type calculi} of second subderivatives and second-order tangents, etc. Furthermore, parabolic regularity married to the {\em metric subregularity constraint qualification} allows us to efficiently deal with {\em constraint systems} in optimization with obtaining {\em computation formulas} for major second-order generalized derivatives, which are instrumental, in particular, for the development and justification of primal-dual numerical algorithms of optimization. Let us also mention in the list of achievements of \cite{mms19} that parabolic regularity leads us to deriving {\em no-gap second-order optimality conditions} in constrained optimization and establishing a {\em quadratic growth} of {\em augmented Lagrangians} that are highly important for subsequent numerical applications. We also refer the reader to \cite{mms19a} for related results in composite models of constrained optimization.\vspace*{0.03in}

In this paper we continue the study and applications of parabolic regularity in the new directions that were not explored in \cite{mms19a,mms19}. Our {\em first goal} is to obtain complete {\em second-order characterizations} of some major stability properties and their {\em robust} versions for perturbed {\em Karush-Kuhn-Tucker} (KKT) systems associated with optimization problems whose constraints are described by parabolically regular sets. We mainly concentrate here on two equivalent (up to taking the inverse mapping) properties known as {\em strong metric subregularity} and {\em isolated calmness} together with their robust counterparts. Using parabolic regularity, we establish novel characterizations of these properties in the general KKT setting under consideration that cover known results for nonlinear programs (NLPs) while contain essentially new ingredients for nonpolyhedral constraint systems. Note that our approach to characterizing isolated calmness of KKT systems in nonpolyhedral constrained optimization, although using some elements similar to the case of polyhedral NLPs, heavily employs advanced machinery of second-order generalized differentiation, which allows us to shed new light on this issue even for classical nonlinear programs. Observe to this end that the recent characterizations of {\em noncritical Lagrange multipliers} obtained in \cite{ms18} for nonpolyhedral systems play a significant role in deriving the aforementioned stability results for perturbed KKT systems and their numerical applications discussed below.\vspace*{0.03in}

The {\em second goal} of this paper is to develop applications of the obtained second-order characterizations of robust stability for KKT systems to the {\em sequential quadratic programming} (SQP) method for solving optimization problems with parabolically regular constraints. In the case of NLPs, the sharpest result for the method was achieved by Bonnans \cite{bo94} who established the {\em superlinear convergence} of the basic SQP method under the uniqueness of Lagrange multipliers and the fulfillment of the classical second-order sufficient condition. We now significantly extend this result to the general class of parabolically regular constrained optimization problems under an appropriate version of the second-order sufficiency with adding a new condition, which automatically holds for NLPs.\vspace*{0.03in}

The rest of the paper is organized as follows. Section~\ref{sect2} presents and discusses the {\em basic concepts} of first-order and second-order variational analysis and generalized differentiation that are systematically used in the subsequent material. In Section~\ref{sect2a} we define the crucial notion of {\em set parabolic regularity}, formulate important results for it taken from \cite{mms19} and needed in this paper, and then derive some new assertions, which are required for our main achievements below.

In Section~\ref{sect3} we first reveal close relationships between {\em noncriticality} of Lagrange multipliers and {\em strong metric subregularity} of the KKT mapping associated with the constrained optimization problem under consideration. It turns out to be instrumental for deriving new second-order characterizations of the isolated calmness property for the solution map to the canonically perturbed KKT system. The latter result is used further in this section to establish {\em robust stability} of the extended KKT system associated with {\em general perturbations} of the original constrained optimization problem. In the case of {\em canonical perturbations}, it gives us the property of solution maps to KKT systems that is labeled as {\em robust isolated calmness} in \cite{dsz}.

Section~\ref{sect4} is devoted to the application of the obtained stability results to establishing local {\em superlinear convergence} of the {\em basic SQP method} to numerically solve optimization problems with {\em parabolically regular constraints}. The stability result under general perturbations developed in Section~\ref{sect3} plays a crucial role in justifying the {\em solvability} of subproblems in the basic SQP method without which the SQP iterations are not well posed. Achieving it and imposing the obtained second-order sufficient condition together with the uniqueness of Lagrange multipliers and the calmness of the multiplier mappings (the latter is automatic for NLPs), we prove superlinear convergence of the basic SQP method and thus properly extend Bonnans' remarkable result to a significantly more general class of nonpolyhedral problems in constrained optimization.

The concluding Section~\ref{sect5} summarizes the main results obtained in the paper and discusses some directions of the future research.

\section{Basic Concepts and Discussions}\sce \label{sect2}

The notation and terminology in this paper are standard in variational analysis and generalized differentiation; see, e.g., \cite{mor18,rw}. For the reader's convenience and notational unification we use as a rule small Greek letters to signify scalar and extended-real-valued (i.e., with values in the extended real line $\oR:=\R\cup\{\infty\}$) functions, small Latin letters for vectors and single-valued mappings, and capital letters for sets and set-valued mappings. For a nonempty subset $\O$ of $\R^n$, the symbols $\ri\O$, $\cl\O$, and $\O^*$ stand for the relative interior, closure, and polar of $\O$,  respectively. We write $x\st{\O}{\to}\ox$ to indicate that $x\to\ox$ with $x\in\O$. The distance between a point $x\in\R^n$ and a set $\O\subset\R^n$ is denoted by ${\rm dist}(x;\O)$. The indicator function of $\O$ is defined by $\dd_\O(x):=0$ for $x\in\O$ and $\dd_\O(x):=\infty$ otherwise. The notation $\B$ stands for the closed unit ball in the space in question, while $\B_r(x):=x+r\B$ mean the closed ball centered at $x$ with radius $r>0$. The symbol $x=o(t)$ with $x\in \R^n$ and $t>0$ tells us that ${\|x\|}/{t}\to 0$ as $t\dn 0$. Let us also mention that the notation $F\colon\R^n\tto\R^m$ indicates the possibility of set values $F(x)\subset\R^m$ (including the empty set $\emp$) of $F$ for some $x\in\R^n$, in contrast to the standard notation $f\colon\R^n\to\R^m$ for single-valued mappings as well as extended-real-valued functions.\vspace*{0.03in}

Now we recall those constructions of generalized differentiation for sets, set-valued mappings, and extended-real-valued functions that are broadly used in what follows. By generalized differential constructions for sets we understand tangent and normal cones, which are closely related to (in fact induce) the corresponding generalized derivatives for mappings and functions.

Given a set $\O\subset\R^n$, the (Bouligand-Severi) {\em tangent cone} $T_\O(\ox)$ to $\O$ at $\ox\in\O$ is defined by
\begin{eqnarray}\label{2.5}
T_\O(\ox):=\big\{w\in\R^n \;\big|\;\exists\,t_k{\downarrow}0,\;\;w_k\to w\;\;\mbox{ as }\;k\to\infty\;\;\mbox{with}\;\;\ox+t_kw_k\in\O\big\}.
\end{eqnarray}
A tangent vector $w\in T_\O(\ox)$ is said to be {\em derivable} if there exists $\xi:[0,\ve]\to \O$ with $\ve>0$, $\xi(0)=\ox$, and $\xi'_+(0)=w$, where $\xi'_+$ signifies the right derivative of $\xi$ at $0$ defined by
$$
\xi'_+(0):=\lim_{t\dn 0}\frac{\xi(t)-\xi(0)}{t}.
$$
The set $\O$ is {\em geometrically derivable} at $\ox$ if every tangent vector $w$ to $\O$ at $\ox$ is derivable.

The (Mordukhovich limiting) {\em normal cone} to $\O$ at $\ox\in\O$ is given by
\begin{eqnarray}\label{2.4}
N_\O(\ox)=\big\{v\in\R^n\;\big|\;\exists\,x_k{\to}\ox,\;\;v_k\to v\;\;\mbox{with}\;\;v_k\in T^*_\O(x_k)\big\},
\end{eqnarray}
where $T^*_\O(x):=\{v\in\R^n\;|\;\la v,w\ra\le 0\;\mbox{ for all }\;w\in T_\O(x)\}$ is the dual cone of \eqref{2.5}. If $\O$ is convex, \eqref{2.5} and \eqref{2.4} reduce, respectively, to the tangent and normal cones of convex analysis. Due to the intrinsic nonconvexity of the normal cone even in simple settings (i.e., for the graph of the function $|x|$ at $(0,0)\in\R^2$), and thus the normal cone \eqref{2.5} is not tangentially generated.

Consider next a set-valued mapping/multifunction $F\colon\R^n\tto\R^m$ with its domain and graph given, respectively, by
$$
\dom F:=\big\{x\in\R^n\;\big|\;F(x)\ne\emp\big\}\quad\mbox{and}\quad \gph F:=\big\{(x,y)\in\R^n\times\R^m\;\big|\;y \in F(x)\big\}.
$$
Then we define the two generalized derivatives generated by the tangent cone \eqref{2.5} and the normal cone \eqref{2.4} to the graph of $F$ in the following way, respectively. The {\em graphical derivative} to $F$ at $(\ox,\oy)\in\gph F$ is
\begin{equation}\label{gder}
DF(\ox,\oy)(u):=\big\{v\in\R^m\;\big|\;(u,v)\in T_{\scriptsize{\gph F}}(\ox,\oy)\big\},\quad u\in\R^n,
\end{equation}
while the {\em coderivative} of $F$ at this point is
\begin{equation}\label{cod}
D^*F(\ox,\oy)(v):=\big\{u\in\R^n \;\big|\;(u,-v)\in N_{\scriptsize{\gph F}}(\ox,\oy)\big\},\quad v\in\R^m.
\end{equation}
Note that the coderivative \eqref{cod} is not dual to the graphical derivative \eqref{gder}, since the values of \eqref{cod} are often nonconvex. Despite (actually due to) its nonconvexity, the coderivative enjoys {\em full calculus} based on the variational/exttemal principles of variational analysis; see \cite{mor06,mor18} for more details. Furthermore, in its terms we get complete pointwise characterizations of fundamental {\em well-posedness} properties of set-valued mappings related to metric regularity, linear openness/covering, and robust Lipschitzian stability known as the Mordukhovich coderivative criteria; see \cite{mor93,mor06,rw}. We present below the one for metric regularity used in this paper. Recall that a mapping $F\colon\R^n\tto\R^m$ is {\em metrically regular} around $(\ox,\oy)\in\gph F$ if there exist neighborhoods $U$ of $\ox$ and $V$ of $\oy$ together with a number $\kappa\ge 0$ ensuring the estimate
\begin{equation}\label{me2}
{\rm dist}\big(x;F^{-1}(y)\big)\le\kappa\,{\rm dist}\big(y;F(x)\big)\quad\mbox{for all}\quad(x,y)\in U\times V.
\end{equation}
The {\em coderivative criterion} for metric regularity says that a closed-graph mapping $F$ is metrically around the point $(\ox,\oy)$ if and only if we have the kernel condition
\begin{equation}\label{cod-cr}
{\rm ker}\,D^*F(\ox,\oy):=\big\{v\in\R^m\;\big|\;0\in D^*F(\ox,\oy)(v)\big\}=\big\{0\big\};
\end{equation}
see \cite[Theorem~3.3]{mor18} for a simplified proof and discussions in finite dimensions.

Along with the usage of \eqref{me2} and \eqref{cod-cr}, we largely employ the following significant modifications and relaxations of \eqref{me2}. The mapping $F$ is said to be {\em metrically subregular} at $(\ox,\oy)$ if the validity of \eqref{me2} is required only when $y=\oy$, i.e., if there exist $\kappa\ge 0$ and a neighborhood $U$ of $\ox$ ensuring the validity of the distance estimate
\begin{equation}\label{metsub}
{\rm dist}\big(x;F^{-1}(\oy)\big)\le\kappa\,{\rm dist}\big(\oy;F(x)\big)\quad\mbox{for all}\quad x\in U.
\end{equation}
Furthermore, $F$ is {\em strongly metrically subregular} at $(\ox,\oy)$ if \eqref{metsub} holds with $F^{-1}(\oy)=\{ \ox\}$  , i.e.,
\begin{equation}\label{str-reg}
\|x-\ox\|\le\kappa\,{\rm dist}\big(\oy;F(x)\big)\quad\mbox{for all}\quad x\in U.
\end{equation}

It has been well realized in variational analysis that the subregularity properties \eqref{metsub} and \eqref{str-reg} are equivalent to the following calmness and isolated calmness counterparts for inverse mappings; see, e.g., \cite{dr,mor18}. Recall that $F\colon\R^n\tto\R^m$ is {\em calm} at $(\ox,\oy)\in\gph F$ if there exist $\ell\ge 0$ and neighborhoods $U$ of $\ox$ and $V$ of $\oy$ such that
\begin{equation}\label{calm}
F(x)\cap V\subset F(\ox)+\ell\|x-\ox\|\B\quad\mbox{for all}\;\;x\in U.
\end{equation}
We say that the mapping $F$ has the {\em isolated calmness} property at $(\ox,\oy)$ if $F(\ox)=\{\oy\}$ in \eqref{calm}, i.e.,
there exist $\ell\ge 0$ and neighborhoods $U$ of $\ox$ and $V$ of $\oy$ for which
\begin{equation}\label{icalm}
F(x)\cap V\subset\big\{\oy\big\}+\ell\|x-\ox\|\B\quad\mbox{for all}\;\;x\in U.
\end{equation}
Although properties \eqref{metsub}--\eqref{icalm} are  much less investigated and applied in the literature in comparison with metric regularity and related two-point properties of multifunctions, there is an increasing number of publications dealing with such one-point formulations; see, e.g., \cite{cdk18,cst,dsz,dr97,dr,g,hms,ho,is14,kk,mms19,mor18,rob81,yy,zn} and the references therein. In this paper we pay the main attention to the properties of  strong metric subregularity and isolated calmness for KKT systems with the subsequent applications to the SQP method.

Given next an extended-real-valued function $\ph\colon\R^n\to\oR$ with the domain $\dom\ph:=\{x\in\R^n\;|\;\ph(x)<\infty\}$ and the epigraph $\epi\ph:=\{(x,\mu)\in\R^{n+1}\;|\;\mu\ge\ph(x)\}$, suppose in what follows that $\ph$ is lower semicontinuous (l.s.c.) around the reference point. The {\em subdifferential} (collection of subgradients) of $\ph$ at $\ox$ is defined by
\begin{equation}\label{sub}
\partial\ph(\ox):=\big\{v\in\R^n\;\big|\;(v,-1)\in N_{\scriptsize{\epi\ph}}\big(\ox,\ph(\ox)\big)\big\}
\end{equation}
via the normal cone \eqref{2.4} to the epigraph of $\ph$ at $(\ox,\ph(\ox))$. If $\ph$ is convex, the subdifferential \eqref{sub} reduces to the classical subgradient set of convex analysis.\vspace*{0.03in}

We proceed further with the second-order constructions of generalized differentiation needed in this paper. Given $\ph\colon\R^n\to\oR$, consider the parametric family of {\em second-order difference quotients} for $\ph$ at $\ox\in\dom\ph$ relative to $\ov\in\R^n$ by
\begin{equation*}
\Delta_t^2\ph(\bar x,\ov)(w):=\dfrac{\ph(\ox+tw)-\ph(\ox)-t\langle\ov,w\rangle}{\sm t^2}\quad\mbox{with}\;\;w\in\R^{n}\;\mbox{ and }\;\;t>0.
\end{equation*}
The {\em second subderivative} of $\ph$ at $\ox$ for $\ov$ and $w$ is defined by
\begin{equation}\label{ssd}
\d^2\ph(\bar x,\ov)(w):=\liminf_{\substack{t\dn 0\\w'\to w}}\Delta_t^2\ph(\ox,\ov)(w').
\end{equation}
Then $\ph$ is said to be {\em twice epi-differentiable} at $\bar x$ for $\ov$ if the sets $\epi\Delta_t^2\ph(\ox,\ov)$ converge (in the sense of the convergence of the corresponding distance functions) to $\epi\d^2\ph(\bar x,\ov)$ as $t\downarrow 0$. If in addition
the second subderivative is a proper function (i.e., does not take the value $-\infty$ and is finite at some point), we say that $\ph$ is {\em properly twice epi-differentiable} at $\bar x$ for $\ov$. By \cite[Proposition~7.2]{rw}, the twice epi-differentiability of $f$ at $\bar x$ for $\ov$  can be understood equivalently as for every $w\in\R^n$ and every sequence $t_k\downarrow 0$, there exists a sequence $w_k\to w$ such that
\begin{equation*}
\Delta_{t_k}^2\ph(\bar x,\ov)(w_k)\to\d^2\ph(\bar x,\ov)(w).
\end{equation*}

\section{Parabolic Regularity and Related Issues}\sce\label{sect2a}

We begin with the formulation of the fundamental parabolic regularity property for sets.

\begin{Definition}[\bf parabolic regularity of sets]\label{par-reg} A set $\O\subset\R^n$ is called  {\sc parabolically regular} at $\ox\in\O$ for a vector $\ov\in\R^n$ if for each $w$ with $\d^2\dd_\O(\bar x,\ov)(w)<\infty$ there exist, among sequences $t_k\dn 0$ and $w_k\to w$ satisfying
\begin{equation*}
\Delta_{t_k}^2\dd_\O(\bar x,\ov)(w_k)\to\d^2\dd_\O(\bar x,\ov)(w)\;\mbox{ as }\;k\to\infty,
\end{equation*}
those with the additional property that
\begin{equation*}
\limsup_{k\to\infty}\frac{\|w_k-w\|}{t_k}<\infty.
\end{equation*}
\end{Definition}

It is shown in \cite{mms19} that this property holds for sets appearing frequently in constrained optimization such as polyhedral convex sets, the second-order/Lorentz/ice-cream cone, and the cone of symmetric positive semidefinite matrices. More broadly, any ${\cal C}^2$-{\em cone reducible set} in the sense of \cite[Definition~3.135]{bs} is always parabolically regular; see \cite[Theorem~6.2]{mms19}. However, the converse statement fails even for simple examples as demonstrated in the following example taken from  \cite{mms19}. Indeed, it is shows in \cite[Example~6.4]{mms19} that for any $\al\in(1,2)$  the set $\O:=\epi g\subset\R^2$ with the function $g\colon\R\to\R$ defined by
\begin{equation*}
g(x):=\left\{\begin{matrix}
0&x\le 0,\\
x^{\alpha}&x\ge 0,
\end{matrix}\right.
\end{equation*}
is not parabolically regular at $\ox=(0,0)\in\O$ and $\ov=(0,-1)$; see the aforementioned example for more detail. One of the  main differences between  the ${\cal C}^2$-cone reducibility and parabolic regularity lies in the fact that the former demands second-order variational information for {\em all tangent vectors} $w\in T_\O(\ox)$, which may not be available. On the contrary, parabolic regularity requires second-order variational information only for tangent vectors that are {\em inside the critical cone} (see \eqref{cri3} for its definition) to the set at the point under consideration. The collections of critical vectors is usually {\em much smaller} than the whole tangent cone. In this example of $\O$, the tangent cone to this set at $\ox$ is $\R\times\R_{+}$, while its critical cone at $\ox$ for $\ov$ is $\R\times\{0\}$. As shown in \cite{mms19}, the second-order variational analysis of sets can be carried out when the vectors from the critical cone enjoy the second-order regularity  defined in Definition~\ref{par-reg}. It is also worth mentioning that,  contrary to the ${\cal C}^2$-cone reducibility, the parabolic regularity enjoys good calculus including intersection rule under reasonable assumptions.

The {\em second-order tangent set} to $\O$ at $\ox\in\O$ for a tangent vector $w\in T_\O(\ox)$ from \eqref{2.5} is
\begin{equation}\label{2tan}
T_\O^2(\ox,w)=\big\{u\in\R^n\;\big|\;\exists\,t_k{\downarrow}0,\;u_k\to u\;\mbox{ as }\;k\to\infty\;\;\mbox{with}\;\;\ox+t_kw+\sm t_k^2u_k\in\O\big\}.
\end{equation}
A set $\O$ is said to be {\em parabolically derivable} at $\ox$ for $w$ if $T_\O^2(\ox,w)\ne\emp$ and for each $u\in T_\O^2(\ox,w)$ there exist $\ve>0$ and an are $\xi:[0,\ve]\to\O$ with $\xi(0)=\ox$, $\xi'_+(0)=w$, and $\xi''_+(0)=u$, where
\begin{equation*}
\xi''_+(0):=\lim_{t\dn 0}\frac{\xi(t)-\xi(0)-t\xi'_+(0)}{\sm t^2}.
\end{equation*}
The latter property is satisfied in rather general settings; see \cite{rw} and \cite{mms19}. Recall further that the {\em critical cone} for a closed set $\O$ at $x$ for $v$ with $(x,v)\in\gph N_\O$ is defined by
\begin{equation}\label{cri3}
 K_\O(x,v)=T_\O(x)\cap\{v\}^\perp.
\end{equation}

Now we are ready to present an important result, which actually summarizes several achievements from \cite{mms19}. To proceed, recall that a set-valued mapping  $F:\R^n\tto\R^m$ is {\em proto-differentiable} at $\ox$ for $\oy$ with $(\ox,\oy)\in\gph F$ if the set $\gph F$ is geometrically derivable at $(\ox,\oy)$. When this condition holds for $F$, we refer to $DF(\ox,\oy)$ as the {\em proto-derivative} of $F$ at $\ox$ for $\oy$.

\begin{Theorem}[\bf proto-differentiability of normal cone mappings]\label{proto} Let $\O\subset\R^n$ be closed and convex set, and let $\ov\in N_{\O}(\ox)$ with $\ox\in\O$. Assume that the set $\O$ is parabolically derivable at $\ox$ for every vector $w\in K_\O(\ox,\ov)$, and that $\O$ is parabolically regular at $\ox$ for $\ov$. Then the following equivalent conditions are satisfied:

{\bf(i)} The indicator function $\dd_\O$ is twice epi-differentiable at $\ox$ for $\ov$.

{\bf(ii)} The normal cone mapping $N_\O$ is proto-differentiable at $\ox$ for $\ov$ with the proto-derivative
\begin{equation}\label{gdpd}
DN_\O(\ox,\ov)(w)=\sub\big[\sm\d^2\dd_\O(\ox,\ov)\big](w)\quad\mbox{for all}\;\;w\in\R^n.
\end{equation}
Furthermore, we have that $\dom\d^2\dd_\O(\ox,\ov)=K_\O(\ox,\ov)$, and that the second subderivative $\d^2\dd_\O(\ox,\ov)$ is a proper convex function on $\R^n$.
\end{Theorem}
{\bf Proof}. The fulfillment of both conditions (i) and (ii) and their equivalence are established in \cite[Theorem~3.6]{mms19} and \cite[Theorem~3.7]{mms19}. The claimed formula for the domain of the second subderivative is taken from \cite[Theorem~3.3]{mms19}. Finally, it follows from \cite[Proposition~13.20]{rw} that the second subderivative $\d^2\dd_\O(\ox,\ov)$ is a proper convex function. $\h$\vspace*{0.05in}

Exploiting the imposed convexity of the set $\O$ allows us to derive the following enhanced representations of the proto-derivative $DN_\O(\ox,\ow)$ from Theorem~\ref{proto} at the origin.

\begin{Corollary}[\bf proto-derivatives of normal cones at the origin]\label{twi3} Let $\O\subset\R^n$ satisfy all the assumptions of Theorem~{\rm\ref{proto}}. Then we have the formulas
\begin{equation}\label{gdpr}
DN_\O(\ox,\ov)(0)=K_\O(\ox,\ov)^*=N_{K_\O(\ox,\ov)}(0).
\end{equation}
\end{Corollary}
{\bf Proof}. The convexity of $\O$, and hence of $K_\O(\ox,\ov)$, allows us to deduce the second equality claimed in \eqref{gdpr}. We proceed now with justifying the first equality therein. It follows from the representation in \eqref{gdpd} that
\begin{equation*}
DN_\O(\ox,\ov)(0)=\sub\big[\sm\d^2\dd_\O(\ox,\ov)\big](0).
\end{equation*}
Pick any $u\in DN_\O(\ox,\ov)(0)$ and thus get $u\in\sub[\frac{1}{2}\d^2\dd_\O(\ox,\ov)](0)$. Since the second subderivative
$\d^2\dd_\O(\ox,\ov)$ is a convex  function with $\d^2\dd_\O(\ox,\ov)(0)=0$, we obtain
\begin{equation*}
\la u,w-0\ra\le\sm\d^2\dd_\O(\ox,\ov)(w)-\sm\d^2\dd_\O(\ox,\ov)(0)=\sm\d^2\dd_\O(\ox,\ov)(w)\quad\mbox{for all}\;\;w\in\R^n.
\end{equation*}
Theorem~\ref{proto} tells us that $\dom\d^2\dd_\O(\ox,\ov)=K_\O(\ox,\ov)$.  To proceed further, pick $t>0$ and $w\in K_\O(\ox,\ov)$ and, with taking into account that the second subderivative $\d^2\dd_\O(\ox,\ov)$ is positive homogeneous of degree $2$, arrive at $\d^2\dd_\O(\ox,\ov)(tw)=t^2\d^2\dd_\O(\ox,\ov)(w)$, which results in
\begin{equation*}
t\la u,w\ra\le\sm\d^2\dd_\O(\ox,\ov)(tw)=\sm t^2\d^2\dd_\O(\ox,\ov)(w).
\end{equation*}
This along with the choice of $w\in\dom\d^2\dd_\O(\ox,\ov)$ implies that $\la u,w\ra\le 0$ and hence gives us $u\in K_\O(\ox,\ov)^*$. This verifies the inclusion ``$\subset$" for the first equality  in \eqref{gdpr}.

Assuming now that $u\in K_\O(\ox,\ov)^*$, we get $\la u,w\ra\le 0$ for all $w\in K_\O(\ox,\ov)$. This leads us to
\begin{equation*}
\la u,w\ra\le 0\le\sm\d^2\dd_\O(\ox,\ov)(w)=\sm\d^2\dd_\O(\ox,\ov)(w)-\sm\d^2\dd_\O(\ox,\ov)(0)\quad\mbox{for all}\;\;w\in K_\O(\ox,\ov).
\end{equation*}
If $w\notin K_\O(\ox,\ov)$, then it follows from Theorem~\ref{proto} that $\d^2\dd_\O(\ox,\ov)(w)=\infty$, which clearly ensures that the above   inequality holds for such $w$. Then using the convexity of second subderivative $\d^2\dd_\O(\ox,\ov)$ and the subdifferential construction of convex analysis, we arrive at $u\in\sub\big[\sm\d^2\dd_\O(\ox,\ov)\big](0)$ and thus deduce from \eqref{gdpd} that $u\in DN_\O(\ox,\ov)(0)$. This verifies the inclusion ``$\supset$" for the first equality in \eqref{gdpr} and hence completes the proof of the corollary. $\h$\vspace*{0.05in}

Next we start considering the {\em constrained optimization} problem given by
\begin{equation}\label{coop}
\min_{x\in\R^n}\;\;\;\ph_0(x)\quad\mbox{subject to }\;\;f(x)\in\Th,
\end{equation}
where $\ph_0\colon\R^n\to\R$ and $f\colon\R^n\to\R^m$ are ${\cal C}^2$-smooth mappings around the reference points, and where $\Th$ is a closed and convex set in $\R^m$. Define its {\em KKT system} associated with \eqref{coop} by
\begin{equation}\label{VS}
\nabla_xL(x,\lm)=\nabla\ph_0(x)+\nabla f(x)^*\lm=0,\quad\;\lm\in N_\Th\big(f(x)\big),
\end{equation}
where $L(x,\lm):=\ph_0(x)+\la f(x),\lm\ra$ is the {\em Lagrangian} associated with \eqref{coop} and $(x,\lm)\in\R^n\times\R^m$.
Given a point $\ox\in\R^n$, denote the set of {\em Lagrange multipliers} at $\ox$ by
\begin{equation}\label{laset}
\Lambda(\ox):=\big\{\lm\in\R^m\;\big|\;\nabla_x L(\ox,\lm)=0,\;\lm\in N_\Th\big(f(\ox)\big)\big\}.
\end{equation}
Having $(\ox,\olm)$ as a solution to the KKT system \eqref{VS} yields $\olm\in\Lambda(\ox)$. It is not hard to see that if  $\olm\in\Lambda(\ox)$, then $\ox$ is a {\em stationary point} of \eqref{VS} in the sense that it satisfies the condition
\begin{equation}\label{stat}
0\in\partial\big(\ph_0+(\dd_\Th\circ g)\big)(\ox),
\end{equation}
where the subdifferential of the generally nonconvex composition in \eqref{stat} is taken from \eqref{sub}.

Now it is time to formulate the {\em standing assumptions} on the closed and convex set in \eqref{coop} that are used in the rest of the paper unless otherwise stated.\\[1.3ex]
{\bf(H1)} For every $\lm\in\Lambda(\ox)$ the set $\Th$ is {\em parabolically derivable} at $f(\ox)$ for all the vectors in the form $\nabla f(\ox)w\subset K_\Th(f(\ox),\lm)$ with some $w\in\R^n$.\\
{\bf(H2)} The set $\Th$ is {\em parabolically regular} at $f(\ox)$ for every $\lambda\in\Lambda(\ox)$.\vspace{0.05in}

As mentioned above, both assumptions (H1) and (H2) hold for important classes of convex sets that naturally and frequently appear in constrained optimization; see, in particular, \cite[Theorem~6.2]{mms19}. Imposing these assumptions on the set $\Th$ opens the door for deriving enhanced second-order optimality conditions for general problems of constrained optimization. The next result is obtained in \cite[Theorem~7.1]{mms19} as a consequence of the comprehensive calculus rules achieved therein for the second subderivatives.

\begin{Theorem}[\bf no-gap second-order optimality conditions]\label{sooc} Let $\ox$ be a feasible solution to problem \eqref{coop} under the validity of the standing assumptions {\rm (H1)} and {\rm (H2)} on $\Th$ with $\ov:=-\nabla\ph_0(\ox)$. Assume further that the mapping $x\mapsto f(x)-\Th$ is metrically subregular at $(\ox,0)$ and consider the set of feasible solutions  to \eqref{coop} given by
\begin{equation}\label{feas}
\O:=\big\{x\in\R^n\;\big|\;f(x)\in\Th\big\}.
\end{equation}
The following second-order optimality conditions for problem \eqref{coop} hold:\\[1ex]
{\bf(i)} If $\ox$ is a local minimum of \eqref{coop}, then the second-order necessary condition
\begin{equation*}
\max_{\lm\in\Lambda(\ox)}\big\{\langle\nabla_{xx}^2L(\bar x,\lm)w,w\rangle+\d^2\dd_\Th\big(f(\bar x),\lm\big)\big(\nabla f(\ox)w\big)\big\}  \ge 0
\end{equation*}
is satisfied for all the critical vectors $w\in K_\O(\ox,\ov)$.\\[1ex]
{\bf(ii)} The fulfillment of the no-gap second-order condition
\begin{equation*}\label{sscc}
\max_{\lm\in\Lambda(\ox)}\big\{\langle\nabla_{xx}^2L(\bar x,\lm)w,w\rangle+\d^2\dd_\Th\big(f(\bar x),\lm\big)\big(\nabla f(\ox)w\big)\big\}>0 \quad\mbox{for all}\quad w\in K_\O(\ox,\ov)\setminus\{0\}
\end{equation*}
amounts to the existence of positive constants $\ell$ and $\ve$ such that the quadratic growth condition
\begin{equation*}\label{quadg2}
\ph(x)\ge\ph(\ox)+\frac{\ell}{2}\|x-\ox\|^2\;\mbox{ for all }\;x\in\B_{\ve}(\ox)\;\mbox{ with }\;\ph:=\ph_0+(\dd_\Th\circ f),
\end{equation*}
which implies that $\ox$ is a strict local minimizer for \eqref{coop}.
\end{Theorem}

Fix $\ox\in\R^n$ and define the {\em multiplier mapping} $M_{\ox}\colon\R^n\times\R^m\tto\R^m$ associated with $\ox$ by
\begin{equation}\label{lagmap}
M_{\ox}(v,w):=\big\{\lm\in\R^m\;\big|\;(v,w)\in G(\ox,\lm)\big\}\quad\mbox{for}\;\;(v,w)\in\R^n\times\R^m,
\end{equation}
which can be viewed as the canonical perturbation $(v,w)$ of the unparameterized {\em KKT mapping} $G\colon\R^n\times\R^m\tto\R^n\times\R^m$ given by
\begin{equation}\label{GKKT}
G(x,\lm):=\left[\begin{array}{c}
\nabla_x L(x,\lm)\\
-f(x)
\end{array}
\right]+\left[\begin{array}{c}
0\\
N_\Th^{-1}(\lm)
\end{array}
\right].
\end{equation}
It is easy to see that $M_{\ox}(0,0)=\Lambda(\ox)$, where $\Lambda(\ox)$ the set of Lagrange multipliers at $\ox$ taken from \eqref{laset}.
We collect some important properties of the multiplier mapping in the next theorem.

\begin{Theorem}[\bf properties of multipliers and qualification conditions]\label{unique} Let $(\ox,\olm)$ be a solution to the KKT system \eqref{VS}. The following properties are equivalent:\\[1ex]
{\bf(i)} The multiplier mapping ${M}_{\ox}$ is calm at $((0,0),\olm)$ and $\Lambda(\ox)=\{\olm\}$.\\[1ex]
{\bf(ii)} The multiplier mapping $M_{\ox}$ enjoys the isolated calmness property at $((0,0),\olm)$.\\[1ex]
{\bf(iii)} We have the dual qualification condition
\begin{equation}\label{gf02}
DN_\Th\big(f(\ox),\olm\big)(0)\cap\ker\nabla f(\ox)^*=\{0\}.
\end{equation}
If in addition the convex set $\Th$ from \eqref{coop} satisfies the standing assumptions {\rm(H1)} and {\rm(H2)}, then the above conditions are equivalent to:\\[1ex]
{\bf(iv)} The strict Robinson constraint qualification
\begin{equation}\label{sqc}
\nabla f(\ox)\R^n+\big[T_\Th\big(f(\ox)\big)\cap\{\olm\}^\bot\big]=\R^m.
\end{equation}
\end{Theorem}
{\bf Proof}. The equivalence between (i), (ii), and (iii) follows from \cite[Theorem~3.1]{ms18}. Furthermore, it is shown in \cite[Proposition~4.3]{ms18} that conditions (iii) and (iv) are equivalent for any ${\cal C}^2$-cone reducible set. The same arguments together with the representation $DN_\Th(f(\ox),\olm)(0)=K_\Th(f(\ox),\olm)^*$ obtained in Corollary~\ref{twi3} under the weaker assumptions (H1) and (H2) can be utilized to verify the claimed equivalence between (iii) and (iv). $\h$\vspace*{0.05in}

We end this section by comparing the dual qualification condition \eqref{gf02} with the {\em basic constraint qualification} for \eqref{coop}, which is also known as the {\em Robinson constraint qualification} for problems of conic programming.

\begin{Proposition}[\bf comparing constraint qualifications]\label{axil2} Let $(\ox,\olm)$ solve the KKT system \eqref{VS}, and let $\Th$ in \eqref{coop} satisfy the standing assumptions {\rm(H1)} and {\rm(H2)}. If the dual qualification condition \eqref{gf02} holds, then the following  constraint qualifications are also fulfilled:\\[1ex]
{\bf(i)} We have the basic constraint qualification
\begin{equation}\label{rcq}
N_\Th\big(f(\ox)\big)\cap\ker\nabla f(\ox)^*=\{0\}.
\end{equation}
{\bf(ii)} For any vector $w\in\R^n$ with $\nabla f(\ox)w\in K_\Th(f(\ox),\olm)$ we have
\begin{equation*}
N_{K_\Th(f(\ox),\olm)}\big(\nabla f(\ox)w\big)\cap\ker\nabla f(\ox)^*=\{0\}.
\end{equation*}
\end{Proposition}
{\bf Proof}. It follows from Theorem~\ref{proto} under the standing assumptions (H1) and (H2) that $\dd_C$ is twice epi-differentiable at $f(\ox)$ for $\olm$. Employing this and \eqref{gdpr} gives us $DN_\Th(f(\ox),\olm)(0)=N_{K_\Th(f(\ox),\olm)}(0)$. Since we always have the inclusion
\begin{equation*}
N_\Th\big(f(\ox)\big)\subset\cl\big[N_\Th\big(f(\ox)\big)+\R\olm\big]=K_\Th\big(f(\ox),\olm\big)^*=N_{K_\Th(f(\ox),\olm)}(0),
\end{equation*}
the dual qualification condition \eqref{gf02} ensures that the basic constraint qualification \eqref{rcq} holds, which therefore verifies assertion (i).

Turning to (ii), pick $w\in\R^n$ with $\nabla f(\ox)w\in K_\Th(f(\ox),\olm)$ and get from Corollary~\ref{twi3} that
\begin{equation*}
N_{K_\Th(f(\ox),\olm)}\big(\nabla f(\ox)w\big)\subset N_{K_\Th(f(\ox),\olm)}(0)=DN_\Th\big(f(\ox),\olm\big)(0).
\end{equation*}
Appealing now to \eqref{gf02} justifies (ii) and thus completes the proof of the proposition. $\h$\vspace*{0.05in}

It is important to mention that, by the coderivative criterion \eqref{cod-cr}, the basic constraint qualification condition \eqref{rcq} amounts to saying that the mapping $x\mapsto f(x)-\Th$ is {\em metrically regular} around $(\ox,0)$. This assumption, which has been conventionally involved in the study of constrained optimization problems of type \eqref{coop}, is significantly more restrictive than the {\em metric subregularity} of the mapping $x\mapsto f(x)-\Th$ at $(\ox,0)$ imposed in Theorem~\ref{sooc} as well as in other results of \cite{mms19} and of this paper.\vspace*{0.05in}

Our subsequent analysis and applications in this paper are conducted under the standing assumptions (H1) and (H2). It follows from Theorem~\ref{unique} that in this setting the primal and dual conditions \eqref{gf02} and \eqref{sqc}, respectively, are equivalent. Having it in mind, we refer to both of these conditions as  the {\em strict Robinson constraint qualification}.

\section{Robust Stability of Perturbed KKT Systems}\sce\label{sect3}

This section is devoted to deriving second-order characterizations of {\em robust stability} of KKT systems by which we precisely mean the equivalent robust versions of strong metric regularity at the reference point of the KKT mapping $G(x,\lm)$ from \eqref{GKKT} and of isolated calmness of its inverse mapping  $S\colon\R^n\times\R^m\tto\R^n\times\R^m$ defined by
\begin{equation}\label{mapS}
S(v,w):=\big\{(x,\lm)\in\R^n\times\R^m\;\big|\;(v,w)\in G(x,\lm)\big\}\quad\mbox{with}\;\;(v,w)\in\R^n\times\R^m,
\end{equation}
The latter mapping can be seen as the {\em solution map} to the {\em canonical perturbation} of the original KKT system \eqref{VS}, i.e., as the KKT system associated with the canonically perturbed constrained optimization problem
\begin{equation}\label{pcoop}
\min_{x\in\R^n}\;\;\;\ph_0(x)-\la v,x\ra\quad\mbox{subject to}\;\;f(x)+w\in\Th.
\end{equation}

It follows from \cite{dr97} that for NLPs the isolated calmness property of \eqref{mapS} is characterized by the simultaneous validity of the classical second-order sufficient optimality condition and the uniqueness of Lagrange multipliers. As further shown in \cite{dsz}, a similar result holds for the case of conic programs in \eqref{coop} with a ${\cal C}^2$-cone reducible set $\Th$, provided that the uniqueness of Lagrange multipliers is replaced by the strict Robinson constraint qualification \eqref{sqc}; see Remark~\ref{rob-calm} for more discussions. In this section we prove that the result of \cite{dr97} can be extended to a general class of constrained optimization problems of type \eqref{coop} with {\em parabolically regular} sets $\Th$ under the additional assumption that the {\em multiplier mapping} \eqref{lagmap} is {\em calm}. Moreover, our results reveal that the combination of the uniqueness of Lagrange multipliers and the calmness of the multiplier mapping amounts to the strict Robinson constraint qualification. It shows that the calmness of the multiplier mapping, being automatically satisfied for NLPs, cannot be dismissed in general. The usage of this property for the study of the robust stability of KKT systems and for the subsequent SQP applications is new not only for parabolically regular problems but also for less general classes of nonpolyhedral programs, while this condition has been recently employed in \cite{ms18} to characterize noncriticality of Lagrange multipliers.

We refer the reader to \cite{ms18} for more discussions on and sufficient conditions for the {\em calmness} of the multiplier mapping \eqref{lagmap}. Besides the NLP setting that corresponds to the {\em polyhedrality} of the set $\Th$ in \eqref{coop}, this calmness holds when the $\Th$ is either second-order/Lorentz cone or the cone of semidefinite symmetric matrices under the validity of the {\em strict complementarity} condition, i.e., if there exists $\lm\in\Lm(\ox)$ with $\lm\in\ri N_\Th(f(\ox))$; see \cite[Theorem~5.10]{ms18}. \vspace*{0.05in}

To achieve our goals, we first establish relationships between the isolated calmness of the solution map $S$ from \eqref{mapS} and the concept of noncriticality of Lagrange multipliers of the KKT systems under consideration. Recall the notions of critical and noncritical multipliers taken from \cite[Definition~3.1]{ms17}. Let $(\ox,\olm)$ be a solution to the KKT system \eqref{VS}. Then the multiplier $\olm\in\Lambda(\ox)$ is said to be {\em critical} for \eqref{VS} if there is a nonzero vector $w\in\R^n$  satisfying the inclusion
\begin{equation}\label{crc}
0\in\nabla_{xx}^2L(\ox,\olm)w+\nabla f(\ox)^*DN_\Th\big(f(\ox),\olm\big)\big(\nabla f(\ox)w\big).
\end{equation}
The multiplier $\olm\in\Lambda(\ox)$ is {\em noncritical} for \eqref{VS} if \eqref{crc} admits only the trivial solution $w=0$. The above definitions are far-going extensions of those introduced by
Izmailov in \cite{iz05} for nonlinear programs with equality constraints, where $\Th=\{0\}\subset\R^m$, and thus the appeal to generalized differentiation is not required. Characterizations of noncriticality are given in \cite[Proposition~1.43]{is14} for NLPs  and in \cite[Theorem~4.1]{ms18} for ${\cal C}^2$-cone reducible constrained problems. Now we come up to the following relationships between noncriticality and isolated calmness.

\begin{Proposition}[\bf characterizations of isolated calmness via noncriticality]\label{uplip} Let $(\ox,\olm)$ be a solution to the KKT system \eqref{VS}.
Consider the following statements:\\[1ex]
{\bf(i)} The Lagrange multiplier $\olm$  is noncritical for \eqref{VS}.\\[1ex]
{\bf(ii)} The solution map $S$ from \eqref{mapS} enjoys the isolated calmness property at $((0,0),(\ox,\olm))$.\\[1.5ex]
Then we have the following relationships:\\[1ex]
{\bf(a)} Implication {\rm(ii)}$\implies${\rm(i)} is always fulfilled.\\[1ex]
{\bf(b)} The converse implication {\rm(i)}$\implies${\rm(ii)} also holds provided that $\Lambda(\ox)=\{\olm\}$ and that the multiplier mapping ${M}_{\ox}$ from \eqref{lagmap} is calm at $((0,0),\olm )$.
\end{Proposition}
{\bf Proof}. It is not hard to deduce from the definition of noncriticality that the Lagrange multiplier $\olm$ is noncritical if and only if we have
\begin{equation}\label{cri2}
\begin{cases}
\nabla_{xx}^2L(\ox,\olm)w+\nabla f(\ox)^*u=0,\\
u\in DN_\Th\big(f(\ox),\olm\big)\big(\nabla f(\ox)w\big)
\end{cases}
\Longrightarrow w=0.
\end{equation}
On the other hand, it follows from the Levy-Rockafellar criterion (see, e.g., \cite[Theorem~4G.1]{dr}) that $S$ enjoys the isolated calmness property at $((0,0),(\ox,\olm))$ if and only if $DS((0,0),(\ox,\olm))(0,0)=\{(0,0)\}$, which amounts to saying that
\begin{equation}\label{isoc}
\begin{cases}
\nabla_{xx}^2L(\ox,\olm)w+\nabla f(\ox)^*u=0,\\
u\in DN_\Th\big(f(\ox),\olm\big)\big(\nabla f(\ox)w\big)
\end{cases}
\Longrightarrow w=0,\;u=0.
\end{equation}
It is clear therefore that assertion (a) is satisfied. Suppose further that the multiplier $\olm$ is noncritical, and so the implication in \eqref{cri2} holds. Then Theorem~\ref{unique} ensures the validity of the dual qualification condition \eqref{gf02} under the assumptions made in (b). Picking a pair $(w,u)$ satisfying the left-hand side relations in \eqref{isoc}, we obtain from \eqref{cri2} that $w=0$. This yields
\begin{equation*}
u\in DN_\Th\big(f(\ox),\olm\big)(0)\cap\ker\nabla f(\ox)^*.
\end{equation*}
Appealing now to \eqref{gf02} confirms that $u=0$, which tells us in turn that \eqref{isoc} holds, and hence we complete the proof of assertion (b) and of the whole proposition. $\h$\vspace*{0.05in}

It is important to notice that the imposed calmness of the multiplier mapping \eqref{lagmap} is essential for the validity of (i)$\implies$(ii) in Proposition~\ref{uplip}. Indeed, \cite[Example~4.8]{ms18} provides a simple semidefinite program for which the calmness of the multiplier mapping $M_{\ox}$ fails and the unique Lagrange multiplier therein is noncritical, but the solution map $S$ does not satisfy the isolated calmness property. On the other hand, the polyhedrality of the convex set $\Th$, along with the noncriticality of the unique Lagrange multiplier, ensures the isolated calmness of $S$. This follows from Proposition~\ref{uplip} due to the automatic fulfillment of the calmness property of $M_{\ox}$ under the imposed polyhedrality as a direct consequence of the classical Hoffman lemma. \vspace*{0.05in}

Now we are ready to derive a major result of this section giving us complete second-order characterizations of the isolated calmness property of the KKT solution map \eqref{mapS} under parabolic regularity of the underlying set $\Th$ in \eqref{coop}.

\begin{Theorem}[\bf second-order characterizations of isolated calmness for KKT systems]\label{isos}
Let $(\ox,\olm)$ be a solution to the KKT system \eqref{VS} under the standing assumptions {\rm(H1)} and {\rm(H2)}. Then the following assertions are equivalent:\\[1ex]
{\bf(i)} The solution map $S$ from \eqref{mapS} enjoys the isolated calmness property at $((0,0),(\ox,\olm))$, and $\ox$ is a local minimizer of \eqref{coop}.\\[1ex]
{\bf(ii)} The second-order sufficient optimality condition
\begin{equation}\label{sosc}
\begin{cases}
\langle\nabla_{xx}^2L(\bar x,\olm)w,w\rangle+\d^2\dd_\Th\big(f(\bar x),\olm\big)\big(\nabla f(\ox)w\big)>0\\
\mbox{for all }\;w\in\R^n\setminus\{0\}\;\;\mbox{with}\;\;\nabla f(\ox)w\in K_\Th\big(f(\ox),\olm\big)
\end{cases}
\end{equation}
holds, the multiplier mapping ${M}_{\ox}$ is calm at $((0,0),\olm)$, and $\Lambda(\ox)=\{\olm\}$.\\[1ex]
{\bf(iii)} The second-order sufficient optimality condition \eqref{sosc} and the strict Robinson constraint qualification \eqref{gf02} are satisfied simultaneously.\\[1ex]
{\bf(iv)} The Lagrange multiplier $\olm$ from \eqref{laset} is noncritical for \eqref{VS}, $\ox$ is a local minimizer of \eqref{coop},
the multiplier mapping ${M}_{\ox}$ is calm at $((0,0),\olm)$, and $\Lambda(\ox)=\{\olm\}$.
\end{Theorem}
{\bf Proof}. We begin with verifying implication (i)$\implies$(iv). Since $S$ enjoys the isolated calmness property at $((0,0),(\ox,\olm))$ and $\Lambda(\ox)$ is convex, we get $\Lambda(\ox)=\{\olm\}$. The calmness of the Lagrange multiplier mapping ${M}_{\ox}$ at $((0,0),\olm)$ is a consequence of the isolated calmness of $S$ at $((0,0),(\ox,\olm))$. Finally, Proposition~\ref{uplip}(a) tells us that $\olm$ is noncritical, and so we arrive at (iv).

Next we prove implication (iv)$\implies$(ii). Assume that (iv) holds and deduce from Theorem~\ref{unique} that \eqref{gf02} is satisfied.
This along with Proposition~\ref{axil2}(i) implies that the basic constraint qualification \eqref{rcq} is also satisfied. As mentioned before, the qualification condition \eqref{rcq} guarantees the metric subregularity of the mapping  $x\mapsto f(x)-\Th$ at $(\ox,0)$ required in Theorem~\ref{sooc}. Employing Theorem~\ref{sooc}(i) with $\Lambda(\ox)=\{\olm\}$ tells us that the second-order necessary optimality condition
\begin{equation}\label{nsosc}
\big\la\nabla^2_{xx}L(\ox,\olm)w,w \big\ra+\d^2\dd_\Th\big(f(\bar x),\olm\big)\big(\nabla f(\ox)w\big)\ge 0\quad\mbox{for all}\;\;w\in K_\O\big(\ox,-\nabla\ph_0(\ox)\big)
\end{equation}
fulfills, where $\O$ is taken from \eqref{feas}. Furthermore, if follows from \eqref{rcq} that
\begin{equation*}
w\in K_\O\big(\ox,-\nabla\ph_0(\ox)\big)\iff\nabla f(\ox)w\in K_\Th\big(f(\ox),\olm\big).
\end{equation*}
To verify the second-order sufficient condition \eqref{sosc}, we need to show that the above inequality is strict for any $w\ne 0$. Arguing by contradiction, suppose that there exists a nonzero vector $\ow\in\R^n$ satisfying the equality
\begin{equation*}
\big\la\nabla^2_{xx}L(\ox,\olm)\ow,\ow\big\ra+\d^2\dd_\Th\big(f(\bar x),\olm\big)\big(\nabla f(\ox)\ow\big)=0\;\;\mbox{with}\;\;\nabla f(\ox)\ow\in K_\Th\big(f(\ox),\olm\big).
\end{equation*}
Consider now the auxiliary optimization problem
\begin{equation}\label{ncp}
\min_{w\in\R^n}\;\;\frac{1}{2}\Big(\big\la\nabla^2_{xx}L(\ox,\olm)w,w\big\ra+\d^2\dd_\Th\big(f(\bar x),\olm\big)\big(\nabla f(\ox)w\big)\Big)
\end{equation}
and deduce from \eqref{nsosc} that $\ow$ is an optimal solution to \eqref{ncp}. It follows from Theorem~\ref{proto} that $\dom\d^2\dd_\Th\big(f(\ox),\bar\lm)=K_\Th\big(f(\ox),\olm\big)$. Combining this with Proposition~\ref{axil2}(ii) confirms that
\begin{equation*}
N_{\scriptsize{\dom \d^2\dd_\Th(f(\ox),\bar\lm)}}\big(\nabla f(\ox)w\big)\cap\ker\nabla f(\ox)^*=\{0\}.
\end{equation*}
Since  the second subderivative $\d^2\dd_\Th\big(f(\ox),\bar\lm)$ is a convex function, the latter condition allows us to obtain the subdifferential chain rule
\begin{equation*}
\sub_w\big(\sm\d^2\dd_\Th\big(f(\bar x),\olm\big)\big(\nabla f(\ox)\cdot)\big)(\ow)=\nabla f(\ox)^*\sub\big[\sm\d^2\dd_\Th\big(f(\ox),\bar\lm)\big]\big(\nabla f(\ox)\ow\big).
\end{equation*}
Remembering that $\ow$ is a minimizer for \eqref{ncp} and employing the subdifferential Fermat rule together the chain rule above and formula \eqref{gdpd} yield
\begin{eqnarray*}
0&\in&\nabla_{xx}^2L(\ox,\olm)\ow+\nabla f(\ox)^*\sub\big[\sm\d^2\dd_\Th\big(f(\ox),\bar\lm\big)\big]\big(\nabla f(\ox)\ow\big)\\
&=&\nabla_{xx}^2L(\ox,\olm)\ow+\nabla f(\ox)^*DN_\Th\big(f(\ox),\bar\lm\big)\big(\nabla f(\ox)\ow\big).
\end{eqnarray*}
Thus we arrive at the inclusion of type \eqref{crc}. Since $\ow\ne 0$, it tells us that $\olm$ is a critical multiplier, a contradiction.
This shows that the inequality in \eqref{nsosc} is strict, which proves the validity of the second-order sufficient condition \eqref{sosc} and hence verifies (ii).

The next step is to check the fulfillment of (ii)$\implies$(i). Assuming (ii), we are going to prove that $\olm$ is noncritical, which amounts to show that implication \eqref{cri2} holds.  To proceed, pick any pair $(w,u)$ satisfying the left-hand side relations in \eqref{cri2}.  This brings us to
\begin{equation}\label{ip01}
\langle\nabla^2_{xx}L(\ox,\olm)w,w\rangle+\la u,\nabla f(\ox)w\ra=0\quad\mbox{and}\quad u\in DN_\Th\big(f(\ox),\olm\big)\big(\nabla f(\ox)w\big).
\end{equation}
The second condition in \eqref{ip01} together with \eqref{gdpd} ensures that $u\in\sub\big[\sm\d^2\dd_\Th\big(f(\ox),\bar\lm\big)\big]\big(\nabla f(\ox)w)$. Thus we obtain the relationships
\begin{equation*}
\nabla f(\ox)w\in\dom\d^2\dd_\Th\big(f(\ox),\bar\lm\big)=K_\Th\big(f(\ox),\olm\big).
\end{equation*}
Recalling again that $\d^2\dd_\Th(f(\ox),\bar\lm)$ is a convex function, it follows from the definition of the subdifferential in convex analysis that
\begin{equation*}
\la u,v-\nabla f(\ox)w\ra\le\sm\d^2\dd_\Th\big(f(\ox),\bar\lm\big)(v)-\sm\d^2\dd_\Th\big(f(\ox),\bar\lm\big)\big(\nabla f(\ox)w\big)
\end{equation*}
whenever $v\in\R^m$. Pick further any $\ve\in(0,1)$ and set $v:=(1\pm\ve)\nabla f(\ox)w$. Since the second subderivative is positive homogeneous of degree $2$, we get
\begin{equation*}
\pm\la u,\nabla f(\ox)w\ra\le\frac{\ve\pm 2}{2}\d^2\dd_\Th(f(\ox),\bar\lm\big)\big(\nabla f(\ox)w\big).
\end{equation*}
Letting $\ve\dn 0$ therein clearly gives us the equality $\d^2\dd_\Th(f(\ox),\bar\lm)(\nabla f(\ox)w)=\la u,\nabla f(\ox)w\ra$.
Combining the latter with \eqref{ip01} implies that
\begin{equation*}
\langle\nabla^2_{xx}L(\ox,\olm)w,w\rangle+\d^2\dd_\Th\big(f(\ox),\bar\lm\big)\big(\nabla f(\ox)w\big)=0,\;\;\nabla f(\ox)w\in K_\Th\big(f(\ox),\olm\big),
\end{equation*}
which results in $w=0$ due to the second-order condition \eqref{sscc}. This verifies that the multiplier $\olm$ is noncritical. Appealing now to Proposition~\ref{uplip}(b) tells us that the solution map $S$ enjoys the isolated calmness property at $((0,0),(\ox,\olm))$. We also conclude from \eqref{sosc}, $\Lm(\ox)=\{\olm\}$, and Theorem~\ref{sooc}(ii) that $\ox$ is a local minimizer of \eqref{coop}, which  proves that (i) holds.

To complete the proof of theorem, it remains to observe that the equivalence between (ii) and (iii) is a direct consequence of Theorem~\ref{unique}. $\h$\vspace*{0.05in}

Let us now discuss relationships of the obtained characterizations of isolated calmness for KKT systems with known results on the subject.

\begin{Remark}[\bf comparison with known results on the KKT isolated calmness]\label{rob-calm} {\rm The equivalence between assertion (i) and (iii) of Theorem~\ref{isos} was obtained in \cite[Theorem~24]{dsz} for the case where ${\cal C}^2$-cone reducible sets $\Th$ by using a different approach. Our new proof for the general case is deeply rooted into the recent developments on parabolic regularity and its important consequences in second-order variational analysis developed in\cite{mms19}. Characterizations (ii) and (iv) of Theorem~\ref{isos} did not appear either in \cite[Theorem~24]{dsz} or any other publication on nonpolyhedral problems. These characterizations highlight significant differences in dealing with KKT systems of type \eqref{VS} where $\Th$ is a nonpolyhedral set. If $\Th$ is polyhedron in \eqref{VS}, the optimization problem \eqref{coop} can be reduced to a nonlinear program. It is known from \cite[Theorem~2.6]{dr97} that in the NLP case the uniqueness of Lagrange multipliers together with \eqref{sosc} amounts to the isolated calmness of the solution mapping $S$ from \eqref{mapS}. As argued in \cite{dsz}, a similar result is not expected for nonpolyhedral constrained optimization while it holds if the uniqueness of multipliers is replaced by the strict Robinson constraint qualification \eqref{sqc}; see \cite[Theorem~24]{dsz}. However, the authors of \cite{dsz} did not address the question to understand why the result with the multiplier uniqueness may fail for nonpolyhedral problems. Theorem~\ref{isos}(ii) answers this question by revealing that the calmness of the multiplier mapping $M_{\ox}$ (which is automatic for NLPs) is essential to derive such an equivalence. As shown in \cite[Example~3.4]{ms18}, the calmness of $M_{\ox}$ may fail for a two-dimensional semidefinite program (and hence for a three-dimensional ice-cream program) even in the case where the set of Lagrange multipliers is a singleton. In that example, the equivalences of Theorem~\ref{isos} also fail.

Finally in this remark, observe that the assumptions of \cite[Theorem~24]{dsz} include the basic constraint qualification \eqref{rcq}, which is not required. Indeed, this condition is always satisfied in Theorem~\ref{isos} because of the imposed assumptions therein}.
\end{Remark}\vspace*{-0.05in}

Now we consider a parameterized version of problem \eqref{coop} together with the corresponding parameterized KKT system and establish crucial {\em well-posedness} and {\em robust stability} results for them with respect to {\em general perturbations}. Given a parameter space $\R^d$ and the mappings $\tph_0\colon\R^n\times\R^d\to\R$ and $\tPh\colon\R^n\times\R^d\to\R^m$, assume that both $\tph_0$ and $\tPh$ are ${\cal C}^2$-smooth with respect to $x\in\R^n$ and $p\in\R^d$. Then we associate with \eqref{coop} its perturbed version
\begin{equation}\label{coop2}
\min_{x\in\R^n}\;\;\;\tph_0(x,p)\quad\mbox{subject to }\;\;\tPh(x,p)\in\Th.
\end{equation}
The next theorem is a nonpolyhedral extension of \cite[Theorem~1.21]{is14} obtained there for NLPs.

\begin{Theorem}[\bf well-posedness and robust stability of solution maps under general perturbations]\label{los} Let $(\ox,\olm)$ be a solution to the KKT system \eqref{VS} under the standing assumptions {\rm(H1)} and {\rm(H2)}. Fixed a nominal parameter $\op\in\R^d$, suppose that that
\begin{equation}\label{assp}
\tph_0(\ox,\op)=\ph_0(\ox), \;\;\nabla_x\tph_0(\ox,\op)=\nabla\ph_0(\ox),\;\;\tPh(\ox,\op)=f(\ox),\;\mbox{ and }\;\nabla_x\tPh(\ox,\op)=\nabla f(\ox).
\end{equation}
Then the following assertions hold:\\[1ex]
{\bf(i)} If $\ox$ is a strict local minimizer of \eqref{coop2} with $p=\op$ and if the basic constraint qualification \eqref{rcq} is satisfied,  then for any $p\in\R^d$ sufficiently close to $\op$ problem \eqref{coop2} admits a local minimizer $x_p$ converging to $\ox$ as $p\to\op$.\\[1ex]
{\bf(ii)} If the second-order sufficient condition \eqref{sosc} and the strict Robinson constraint qualification \eqref{gf02} are satisfied,
and if
\begin{equation}\label{assp2}
\nabla^2_{xx}\tph_0(\ox,\op)+\nabla^2_{xx}\la\olm,\tPh\ra(\ox,\op)=\nabla^2_{xx}L(\ox,\olm),
\end{equation}
then there exist constants $\ve>0$ and $\ell\ge 0$ such that
\begin{equation}\label{isop}
\emp\ne\Upsilon(p)\cap\B_\ve(\ox,\olm)\subset\{(\ox,\olm)\}+\ell\|p-\op\|\quad\mbox{for all}\;\;p\in\B_\ve(\op),
\end{equation}
where $\Upsilon:\R^d\tto\R^n\times\R^m$ is the solution map to the KKT system of \eqref{coop2} defined by
\begin{equation*}
\Upsilon(p):=\Big\{(x,\lm)\in\R^n\times\R^m\;\Big|\,
\left[\begin{array}{c}
0\\
0
\end{array}
\right]\in\left[\begin{array}{c}
\nabla_x\tph_0(x,p)+\nabla_x\tPh(x,p)^*\lm\\
-\tPh(x,p)
\end{array}
\right]+\left[\begin{array}{c}
0\\
N_\Th^{-1}(\lm)
\end{array}
\right]
\Big\}.
\end{equation*}
\end{Theorem}
{\bf Proof}. Since $\ox$ is a strict local minimum of \eqref{coop2} as $p=\op$, we get a number $\ve>0$ such that $\ox$ is the unique minimizer for the problem
\begin{equation}\label{pro0}
\min_{x\in\R^n}\;\;\;\tph_0(x,\op)\quad\mbox{subject to}\;\;\tPh(x,\op)\in\Th,\;x\in\B_\ve(\ox).
\end{equation}
Consider the auxiliary constrained optimization problem
\begin{equation}\label{pro1}
\min_{x\in\R^n}\;\;\;\tph_0(x,p)\quad\mbox{subject to}\;\;\tPh(x,p)\in\Th,\;x\in\B_\ve(\ox).
\end{equation}
It follows from \eqref{rcq}  and \eqref{assp} that
\begin{equation}\label{rcqp}
N_\Th \big(\tPh(\ox,\op)\big)\cap\ker\nabla_x\tPh(\ox,\op)^*=\{0\},
\end{equation}
which is equivalent to the Robinson constraint qualification for the constraint system $\tPh(x,p)\in\Th$ at $(\ox,\op)$, i.e., to the condition
\begin{equation*}
0\in{\rm int}\big(\tPh(\ox,\op)+\nabla_x\tPh(\ox,\op)\R^n-\Th\big).
\end{equation*}
Thus it follows from \cite[Theorem~2.86]{bs} that for all $(x,p)$ sufficiently close to $(\ox,\op)$ we have
\begin{equation}\label{pro2}
{\rm dist}\big(x;\Gamma(p)\big)=O\big({\rm dist}(\tPh(x;p);\Th)\big),
\end{equation}
where $\Gamma(p):=\{x\in\R^n\;|\;\tPh(x,p)\in\Th\}$. This guarantees that whenever $p$ is sufficiently close to $\op$ the feasible region of \eqref{pro1}, namely $\Gamma(p)\cap\B_\ve(\ox)$, is nonempty for all small positive numbers $\ve$. Appealing to the classical Weierstrass theorem, we deduce that \eqref{pro1} attains a (global) minimizer $\ox_p$ for all $p$ sufficiently close to $\op$.

Now we claim that $\ox_p\to\ox$ as $p\to\op$. Indeed, the failure of it gives us a sequence $p_k\to\op$ for which $\ox_{p_k}$, an optimal solution of \eqref{pro1} for $p=p_k$, does not converge to $\ox$ as $k\to \infty$. Since $\ox_{p_k}\in\B_\ve(\ox)$, there exists a subsequence of $\ox_{p_k}$ that converges to some $\tilde{x}\in\B_\ve(\ox)$ with $\tilde{x}\ne\ox$. Without relabeling, assume with no loss of generality that $\ox_{p_k}\to\tilde x$ for all $k\to\infty$ and then find by \eqref{pro2} such $\tilde{x}_{p_k}\in\Gamma(p_k)$ that
\begin{equation*}
\|\tilde{x}_{p_k}-\ox\|=O\big({\rm dist}(\tPh(\ox,p_k);\Th)\big)=O(\|p_k-\op\|).
\end{equation*}
Remembering that $p_k\to\op$ as $k\to\infty$, it is no harm to suppose that $\tilde{x}_{p_k}\in\B_\ve(\ox)$ for all $k\in\N$. Since $\ox_{p_k}$ is an optimal solution of \eqref{pro1} for $p=p_k$, we have $\tph_0({\ox}_{p_k},p_k)\le\tph_0(\tilde{x}_{p_k},p_k)$. Passing to the limit as $k\to\infty$ brings us to $\tph_0(\tilde{x},\op)\le\tph_0(\ox,\op)$, which contradicts the fact that $\ox$ is the unique minimizer of problem \eqref{pro0}. Thus we are done with (i).

Turing next to (ii), pick any $\ve>0$ such that (i) holds whenever $p\in\B_\ve(\op)$. We claim that the set of Lagrange multipliers for \eqref{coop2} at $(\ox,\op)$ is the singleton $\{\olm\}$, which follows from the combination of Theorem~\ref{unique} with conditions \eqref{gf02} and \eqref{assp}. This allows us to deduce from \eqref{assp2}, the second-order sufficient condition \eqref{sosc}, and Theorem~\ref{sooc}(ii) that $\ox$ is a strict local minimizer of \eqref{coop2} for $p=\op$. Employing now assertion (i) of this theorem,  we conclude  that for any $p$ sufficiently close to $\op$ problem \eqref{coop2} admits a local minimizer $x_p$ with $x_p\to\ox$ as $p\to\op$. It follows from Proposition~\ref{axil2}(i) combined with conditions \eqref{gf02} and \eqref{assp} that the basic constraint qualification \eqref{rcqp} for problem \eqref{pro0} is satisfied. Shrinking the closed ball $\B_\ve(\op)$ if necessary, we assume without loss of generality that the basic constraint qualification \eqref{rcqp} holds for the perturbed problem \eqref{pro1} at $(x_p,p)$ for all $p\in\B_\ve(\op)$. This ensures the existence of a Lagrange multiplier $\lm_p$ associated with the local minimizer $x_p$ for \eqref{coop2} whenever $p\in\B_\ve(\op)$. Shrinking again the ball $\B_\ve(\op)$, we conclude that the Lagrange multipliers $\lm_p$ are uniformly bounded for all $p\in\B_\ve(\op)$. Remembering that the set of Lagrange multipliers for \eqref{coop2} at $(\ox,\op)$ is the singleton $\{\olm\}$ and using it together with the boundedness of $\{\lm_p\}$ for $p\in\B_\ve(\op)$ yield the convergence $\lm_p\to\olm$ as $p\to\op$. Since $(x_p,\lm_p)\in\Upsilon(p)$ and $(x_p,\lm_p)\to(\ox,\olm)$ as $p\to\op$, by decreasing $\ve>0$ if necessary we get $\Upsilon(p)\cap\B_\ve(\ox,\olm)\ne\emp$ for all $p\in\B_\ve(\op)$.

To justify finally the inclusion in \eqref{isop}, observe first by  \eqref{assp} that $(\ox,\olm)\in\Upsilon(\op)$.
Employing  the latter together with \eqref{assp2} and \cite[Corollary~4E.3]{dr} tells us that the mapping $\Upsilon$ enjoys the isolated calmness property at $(\op,(\ox,\olm))$ provided that the implication
\begin{equation}\label{lr}
\begin{cases}
\nabla_{xx}^2L(\ox,\olm)w+\nabla f(\ox)^*u=0,\\
u\in DN_\Th\big(f(\ox),\olm\big)\big(\nabla f(\ox)w\big)
\end{cases}
\Longrightarrow w=0,\;u=0
\end{equation}
holds. To verify \eqref{lr}, we proceed similarly to the proof of implication (ii)$\implies$(i) in Theorem~\ref{isos} and thus complete the proof of the theorem.$\h$\vspace*{0.05in}

The following remark is important to understand the meaning of the robust stability property established in Theorem~\ref{los}(ii) and its comparison with known results in this direction.

\begin{Remark}[\bf robust stability of KKT systems with respect to general vs.\ canonical perturbations]\label{gen-can}
{\rm The properties established in \eqref{isop} of Theorem~\ref{los}(ii) were investigated before for nonlinear programming problems in \cite[Theorem~2.6]{dr97}. Then they were extended in \cite[Theorem~24]{dsz} to ${\cal C}^2$-cone reducible constrained optimization problems  in the case where the parameterized problem \eqref{coop2} has the {\em canonically perturbed} form \eqref{pcoop}  meaning that
\begin{equation*}
\tph_0(x,p)=\ph_0(x)-\la v,x\ra\;\mbox{ and }\;\tPh(x,p)=f(x)+w\;\mbox{ with }\;p:=(v,w)\in\R^n\times\R^m.
\end{equation*}
The properties in \eqref{isop} of Theorem~\ref{los}(ii) were called in \cite{dsz} the {\em robust isolated calmness}. We can see that the difference between \eqref{isop} and the isolated calmness of the solution map $\Upsilon$ is the additional highly important assertion that $\Upsilon(p)\cap\B_\ve(\ox,\olm)\ne\emp$ for all $p$ sufficiently close to $\op$. In general, the isolated calmness and its robust counterpart \eqref{isop} are {\em not equivalent} as shown in \cite[Example~6.4]{mor} for a three-dimensional second-order cone program.

On the other hand, it follows from the previous Theorem~\ref{isos} that the assumptions imposed in Theorem~\ref{los}(ii) yields the isolated calmness property at $((0,0),(\ox,\olm))$ of the solution map \eqref{mapS} for the KKT system associated with the canonically perturbed form \eqref{pcoop}. Combining it with Theorem~\ref{los}(ii) tells us that for {\em canonical perturbations} the isolated calmness and robust isolated properties of KKT solution maps are {\em equivalent}. Our main reason to consider {\em general perturbations} as in \eqref{coop2} and to justify the robust isolated calmness for the associated KKT systems is that this allows us to verify the {\em solvability} of subproblems in the basic SQP method for parabolically regular constrained optimization problems as developed in the next section. It does not seem to be possible while considering only canonical perturbations.

Note to this end that the assumptions   \eqref{assp} and \eqref{assp2} imposed in Theorem~\ref{los} were only
imposed to avoid adjusting \eqref{rcq},  \eqref{sosc}, and  \eqref{gf02} for the perturbed problem \eqref{coop2} and are
automatically satisfied in our application to the aforementioned solvability of the basic SQP method for the original problem \eqref{coop} that is addressed below in Theorem~\ref{solva}.}
 \end{Remark}

\section{Superlinear Convergence of the Basic SQP Method}\sce\label{sect4}

This section is devoted to establish the {\em primal-dual superlinear convergence} of the basic SQP method for the general class of parabolically regular constrained optimization problems of type \eqref{coop}. Various versions of the SQP method have been among the most effective numerical algorithms in constrained optimization. The principal idea of SQP algorithms is to solve  a sequence of quadratic
approximations, called subproblems, whose optimal solutions converge under appropriate assumptions to an optimal solution of the original problem.

Given the current iterate $x_k$, the {\em generic SQP subproblems} for \eqref{coop} are formulated as
\begin{equation}\label{subs1}
\begin{cases}
\disp{\min_{x\in\R^n}}&\ph_0(x_k)+\la\nabla \ph_0(x_k),x-x_k\ra+\sm\la H_k(x-x_k),x-x_k\ra\\
\mbox{subject to}&f(x_k)+\nabla f(x_k)(x-x_k)\in\Th,
\end{cases}
\end{equation}
where $H_k\in\R^{n\times n}$ for all $k=0,1,\ldots$. The KKT system associated with each subproblem \eqref{subs1} is described by the generalized equation
\begin{equation}\label{kktsub}
\left[\begin{array}{c}
0\\
0
\end{array}
\right]\in
\left[\begin{array}{c}
\nabla_{x}L(x_k,\lm_k)\\
-f(x_k)
\end{array}
\right]+\left[\begin{array}{cc}
H_k&\nabla f(x_k)^*\\
-\nabla f(x_k)&0
\end{array}
\right]
\left[\begin{array}{c}
x-x_k\\
\lm-\lm_k
\end{array}
\right]+\left[\begin{array}{c}
0\\
N_\Th^{-1}(\lm)
\end{array}
\right].
\end{equation}
The {\em generic SQP method} for the original problem \eqref{coop} is therefore as follows:
\begin{Algorithm} [\bf generic SQP method]\label{sqpal} {\rm Choose $(x_k,\lm_k)\in\R^n\times\R^m$ and set $k=0$.\\[1ex]
{\bf 1.} If $(x_k,\lm_k)$ satisfies the KKT system \eqref{VS}, then stop.\\[1ex]
{\bf 2.} Choose some matrix $H_k\in\R^{n\times n}$ and compute the next iterate $(x_{k+1},\lm_{k+1})$ as a solution to the generalized equation \eqref{kktsub}.\\[1ex]
{\bf 3.} Increase $k$ by $1$ and then go back to Step~{\rm 1}.}
\end{Algorithm}
In {\em quasi-Newton SQP} algorithms, the matrix $H_k$ is chosen as a quasi-Newton approximation of the Hessian $\nabla_{xx}^2L(x_k,\lm_k)$.
An efficient way to construct $H_k$ is the {\em BFGS} (Broyden--Fletcher--Goldfarb--Shanno) method; see \cite{fp,is14,kk} for more details on this and related algorithms. Instructive discussion and suggestions on how to construct quasi-Newton SQP approximations can be found in \cite{bgls}. Following \cite{is14}, we refer to the {\em basic SQP method} if $H_k$ are chosen as
\begin{equation}\label{hk}
H_k= \nabla_{xx}^2L(x_k,\lm_k)\;\mbox{ for all }\;k=0,1,\ldots.
\end{equation}

It has been well recognized, the basic SQP method can be viewed as a natural extension of the {\em Newton method} that is implemented for {\em generalized equations} rather than equations.
Indeed, we can equivalently express the KKT system \eqref{VS} as the generalized equation
\begin{equation}\label{GKKT2}
\left[\begin{array}{c}
0\\
0
\end{array}
\right]\in\left[\begin{array}{c}
\nabla_x L(x,\lm)\\
-f(x)
\end{array}
\right]+\left[\begin{array}{c}
0\\
N_\Th^{-1}(\lm)
\end{array}
\right].
\end{equation}
Employing the Newton method for this generalized equation leads us to the aforementioned basic SQP method for the constrained problem \eqref{coop}; see \cite[Section~3.1]{is14} for more details.

The Newton method for generalized equations has been investigated extensively since the late 1970s; see, e.g., \cite{dr,fp,is14,kk} and the references therein. Josephy's observation in \cite{jo} significantly advanced the topic by showing that {\em strong regularity} in the sense of Robinson \cite{rob80} ensures the superlinear convergence of the Newton method for variational inequalities. As proved in \cite[Theorem~4.1]{rob80} for NLPs, strong regularity is guaranteed under a stronger form of the second-order sufficient condition \eqref{sosc} together with the linear independence constraint qualification. The next improvement in the setting of NLPs was achieved by Bonnans in \cite{bo94} by showing that the basic SQP method converges superlinearly if the second-order sufficient condition \eqref{sosc} with a {\em polyhedral} set $\O$ therein is satisfied and the set of Lagrange multipliers is a {\em singleton}.

In this section we are going to extend Bonnans' result to nonpolyhedral constrained optimization problems \eqref{coop} in the general case of {\em parabolically regular} sets $\Th$. To do it, we intend to employ the results from the previous section, which tell us that an extra {\em calmness} assumption on the multiplier mapping \eqref{lagmap} is needed. To achieve our goal, we utilize---besides Theorems~\ref{isos} and \ref{los} obtained above---the result of Theorem~3.2 from the book by Izmailov and Solodov \cite{is14} in which superlinear convergence of the Newton method for abstract generalized equations was established under the following two assumptions. The first one, called ``semistability" in \cite{is14}, reduces for \eqref{GKKT2} to the isolated calmness of the solution map $S$ from \eqref{mapS}. The second assumption of \cite[Theorem~3.2]{is14}, called ``hemistability" therein, is defined in \cite[Definition~3.1]{is14} for arbitrary generalized equations of the type
\begin{equation}\label{ge}
0\in g(x)+F(x),\quad x\in\R^n,
\end{equation}
with $g:\R^n\to\R^m$ and $F:\R^n\tto\R^m$. A feasible solution $\ox$ to \eqref{ge} is said to be {\em hemistable} if for any $x$ close to $\ox$ the perturbed generalized equation
\begin{equation*}
0\in g(x)+\nabla g(x)\eta+F(x+\eta)
\end{equation*}
has a solution $\eta_x$ such that $\eta_x\to 0$ as $x\to\ox$. We show below that the isolated calmness of the solution map \eqref{mapS} at $((0,0),(\ox,\olm))$
when $\ox$ is a local minimum of \eqref{coop} yields the hemistability of the solution pair $(\ox,\olm)$ for \eqref{GKKT2}. It is important to notice that the isolated calmness alone can not ensure the hemistability of generalized equations; see \cite[Example~3.3]{is14}. \vspace*{0.05in}

The first major theorem of this section justifies the crucial {\em well-posedness} of subproblems \eqref{subs1} in the basic SQP case \eqref{hk} under parabolic regularity, in the sense of the existence of optimal solutions to \eqref{subs1} (i.e., {\em solvability} of these subproblems) and {\em stability} of the corresponding KKT systems with respect to general parameter perturbations. The proofs of these results are based on the robust stability characterizations established in Section~\ref{sect3}, which are in turn employ recent developments on parabolic regularity \cite{mms19} and noncriticality of Lagrange multipliers \cite{ms18}.

\begin{Theorem}[\bf solvability and stability of subproblems in the basic SQP method]\label{solva} Let $(\ox,\olm)$ be a solution to the KKT system \eqref{VS} under the standing assumptions {\rm(H1)} and {\rm(H2)}. Suppose in addition that the second-order sufficient condition \eqref{sosc} and the strict Robinson constraint qualification \eqref{gf02} are satisfied. Then there is $\ve>0$ such that for all $(u,\mu)\in\B_\ve(\ox,\olm)$ the following assertions hold:\\[1ex]
{\bf(i)} The perturbed optimization problem
\begin{equation}\label{subs}
\begin{cases}
\disp{\min_{x\in\R^n}}&\ph_0(u)+\la\nabla\ph_0(u),x-u\ra+\sm\la\nabla_{xx}^2L(u,\mu)(x-u),x-u\ra\\
{\rm{subject \,to}}&f(u)+\nabla f(u)(x-u)\in\Th
\end{cases}
\end{equation}
always admits a local optimal solution.\\[1ex]
{\bf(ii)} The KKT system of \eqref{subs}, which can be formulated as
\begin{equation}\label{kktsub4}
\left[\begin{array}{c}
0\\0
\end{array}
\right]\in\left[\begin{array}{c}
\nabla_{x}L(u,\mu)\\
-f(u)
\end{array}
\right]+\left[\begin{array}{cc}
\nabla_{xx}^2L(u,\mu)&\nabla f(u)^*\\
-\nabla f(u)&0
\end{array}
\right]
\left[\begin{array}{c}
x-u\\
\lm-\mu
\end{array}
\right]
+\left[\begin{array}{c}
0\\
N_\Th^{-1}(\lm)
\end{array}
\right],
\end{equation}
has a solution $(x,\lm)$ that converges to $(\ox,\olm)$ as $(u,\mu)\to(\ox,\olm)$.
\end{Theorem}
{\bf Proof}. Set $p:=(u,\mu)$ with $\op:=(\ox,\olm)$, denote
\begin{equation*}
\ph(x,p):=\ph_0(u)+\la\nabla\ph_0(u, x-u\ra+\sm\la\nabla_{xx}^2L(u,\mu)(x-u),x-u\ra,\;\;g(x,p):=f(u)+\nabla f(u)(x-u),
\end{equation*}
and view \eqref{subs} as the parameterized optimization problem
\begin{equation}\label{coop3}
\min_{x\in\R^n}\;\;\;\ph(x,p)\quad\mbox{subject to}\;\;g(x,p)\in\Th
\end{equation}
with respect to the {\em general} parameter perturbations $p$. We claim now that $x:=\ox$ is a {\em strict local minimizer} for \eqref{coop3} associated with $\op$. To this end, observe first that the KKT system of \eqref{coop3} associated with $\op$ has the representation
\begin{equation}\label{kktsub2}
\left[\begin{array}{c}
0\\
0
\end{array}
\right]\in
\left[\begin{array}{c}
\nabla\ph(\ox)+\nabla_{xx}^2L(\ox,\olm)(x-\ox)+\nabla f(\ox)^*\lm\\
-f(\ox)-\nabla(\ox)(x-\ox)
\end{array}
\right]+\left[\begin{array}{c}
0\\
N_\Th^{-1}(\lm)
\end{array}
\right].
\end{equation}
We clearly have that$(\ox,\olm)$ is a solution to \eqref{kktsub2}, which hence implies that $\ox$ is a stationary point for \eqref{coop3} associated with $\op$ and that $\olm$ is a Lagrange multiplier associated with $\ox$ for the latter problem. Furthermore, it is not hard to observe from the KKT system \eqref{kktsub2} with $p=\op$ that the set of Lagrange multipliers for the latter problem at $(\ox,\op)$ agrees with that for problem \eqref{coop} at $\ox$. Combining these observations tells us that the set of Lagrange multipliers for the parameterized problem \eqref{coop3} at $(\ox,\op)$ is the singleton $\{\olm\}$. Since
\begin{equation*}
\nabla^2_{xx}\ph(\ox,\op)+\nabla^2_{xx}\la\olm,g\ra(\ox,\op)=\nabla^2_{xx}L(\ox,\olm)\quad\mbox{and}\quad\nabla_x g(\ox,\op)=\nabla f(\ox),
\end{equation*}
and since the second-order sufficient condition \eqref{sosc} holds for \eqref{coop} at $(\ox,\olm)$, we get
\begin{eqnarray*}
&&\big\la\big[\nabla^2_{xx}\ph(\ox,\op)+\nabla^2_{xx}\la\olm,g\ra(\ox,\op)\big]w,w\big\ra+\d^2\dd_\Th\big(f(\ox,\op),\olm\big)\big(\nabla_x g(\ox,\op)w\big)\\
&&=\big\la\nabla^2_{xx}L(\ox,\olm)w,w\big\ra+\d^2\dd_\Th\big(f(\bar x),\olm\big)\big(\nabla f(\ox)w\big)>0\\
&&\mbox{for all}\;w \in\R^n\setminus\{0\}\;\;\mbox{with}\;\;\nabla_x g(\ox,\op)w=\nabla f(\ox)w\in K_\Th\big(f(\ox),\olm\big).
\end{eqnarray*}
This verifies therefore that the second-order sufficient condition of type \eqref{sosc} for problem \eqref{coop3} is satisfied at $((\ox,\op),\olm)$. Employing further Theorem~\ref{sooc}(ii) ensures that $\ox$ a strict local minimizer for \eqref{coop3}
as $p=\op$. Then it follows from Theorem~\ref{los} that there exists $\ve>0$ such that for all $(u,\mu)=p\in\B_\ve(\op)$ the parameterized problem \eqref{coop3} admits a local minimizer $x_p$ with $x_p\to\ox$ as $p\to\op$. This concludes the proof of assertion (i).

Turing next to (ii), Turing now to (ii), observed first that the KKT system of \eqref{subs} is the same is the KKT system of the parametric problem \eqref{coop3}. By  Theorem~\ref{los}(ii), for any $p=(u,\mu)$ close to  $\op:=(\ox,\olm)$ the generalized equation \eqref{kktsub4} admits a solution $(x_p,\lm_p)$. Moreover, \eqref{isop} ensures that $(x_p,\lm_p)\to(\ox,\olm)$ as $p\to\op$,
 which  completes the proof of (ii). $\h$\vspace*{0.05in}

To the best of our knowledge, the solvability of subproblems in the basic SQP method for NLPs was first established by Robinson \cite[Theorem~3.1]{rob76} who assumed, in addition to the second-order sufficient condition, the validity of the linear independence constraint qualification together with the strict complementarity condition. This result was improved by Bonnans \cite[Proposition~6.3]{bo94} for NLPs satisfying the second-order sufficient condition \eqref{sosc} and having a unique Lagrange multiplier. Theorem~\ref{solva} extends Bonnans' result to a general class of nonpolyhedral problems with parabolically regular constraint sets.\vspace*{0.05in}

We are finally in a position to derive the superlinear convergence of the basic SQP method in our general setting. To be precise, the aimed convergence is understood as the {\em superlinear} convergence of the primal-dual iterates $(x_k,\lm_k)$ to a given KKT solution $(\ox,\olm)$ meaning that
\begin{equation*}
\|(x_{k+1},\lm_{k+1})-(\ox,\olm)\|=o\big(\|(x_k,\lm_k)-(\ox,\olm)\|\big)\;\mbox{ as }\;k\to\infty.
\end{equation*}

\begin{Theorem}[\bf superlinear convergence of the basic SQP method under parabolic regularity]\label{supsqp} Let $(\ox,\olm)$ be a solution to the KKT system \eqref{VS} under the validity of the standing assumptions {\rm(H1)} and {\rm(H2)}. Suppose in addition that  the second-order sufficient condition \eqref{sosc} is satisfied at $(\ox,\olm)$, that the multiplier mapping ${M}_{\ox}$ is calm at $((0,0),\olm)$, and that $\Lambda(\ox)=\{\olm\}$. Then there exists a positive constant $\dd$ such that  for any starting point $(x_0,\lm_0)\in\R^n\times\R^m$ sufficiently close to $(\ox,\olm)$, Algorithm~{\rm\ref{sqpal}} with $H_k$ taken from \eqref{hk} generates an iterative sequence $\{(x_k,\lm_k)\}\subset\R^n\times\R^m$ satisfying
\begin{equation}\label{lc1}
\|(x_{k+1}-x_k,\lm_{k+1}-\lm_k)\|\le\dd
\end{equation}
that converges to $(\ox,\olm)$, and the rate of convergence is superlinear.
\end{Theorem}
{\bf Proof}. Due to our results achieved above, we can fit into the framework of \cite[Theorem~3.2]{is14} for the KKT system \eqref{VS}, which can equivalently written as the generalized equation \eqref{GKKT2}. Indeed, the semistability and hemistability properties of $(\ox,\olm)$ imposed in \cite[Theorem~3.2]{is14} are the isolated calmness of the solution map $S$ and the KKT well-posedness formulated in Theorem~\ref{solva}(ii), respectively. Then both of these requirements are satisfied by the results of Theorems~\ref{isos} and \ref{solva}. It allows us to deduce the claimed superlinear convergence of primal-dual iterates from the abstract result of \cite[Theorem~3.2]{is14} for arbitrary generalized equations. $\h$ \vspace*{0.05in}

Let us now compare Theorem~\ref{supsqp} with known results in the literature.

\begin{Remark}[\bf comparison with known results on superlinear convergence of the basic SQP method]\label{comp-sqp} {\rm We address here the following two issues:

{\bf(i)} Theorem~\ref{supsqp} extends the sharpest local superlinear convergence result obtained in \cite[Theorem~5.1]{bo94} for NLPs to a large class of parabolically regular constrained optimization problems. The additional calmness requirement for the multiplier mapping \eqref{lagmap} automatically holds not only for polyhedral problems, but also in the essentially more general settings; see the  discussion above in front of Proposition~\ref{uplip}. Note also that the results on SQP superlinear convergence for the particular class of second-order cone programs were given in \cite{kf,wz} under the strong regularity of the KKT system \eqref{VS} that is a significantly more restrictive  assumption in comparison with those imposed in Theorem~\ref{supsqp}.

{\bf(ii)} As pointed out by a referee, the localization condition \eqref{lc1} cannot be dropped from Theorem~\ref{supsqp} even for the usual NLP setting. The reader is referred to \cite[Examples~5.1 and 5.2]{is15} for more discussions on this localization condition.

{\bf(iii)} The superlinear convergence of a sequence of the SQP iterates $(x_k,\lm_k)$, while {\em not the existence} of such a sequence, can be derived from \cite[Theorem~6.4]{cdk18} in which the superlinear convergence of the Newton method for the generalized equation \eqref{ge}
was obtained when the mapping $g+F$ is strongly metrically subregular. Adopting the latter result to our framework of the generalized equation \eqref{GKKT2}, the assumed strong metric subregularity is equivalent to the isolated calmness of the solution map $S$ from \eqref{mapS}.
Furthermore, by Theorem~\ref{isos} we can equivalently translate the imposed assumptions in Theorem~\ref{supsqp} to the requirements that the solution map $S$ enjoys the isolated calmness property and the point $\ox$ is a local minimizer for \eqref{coop}. This indicates that, in contrast to \cite{cdk18}, we impose the additional assumption on local optimality of $\ox$, and so one might ask whether it is possible to drop this condition with no harm. However, it happens not to be correct since \cite[Theorem~6.4]{cdk18} {\em assumes} the existence of iterates $(x_k,\lm_k)$ (i.e., {\em solvability} of the subproblems) that always stay in a neighborhood of $(\ox,\olm)$. Below we present an example of an NLP problem considered in \cite[Example~3.3]{is14}, which shows that the isolated calmness of the solution map $S$ alone is not enough to guarantee the solvability of the subproblems in our setting. Thus the {\em local optimality} of $\ox$ turns out to be an {\em essential} requirement to construct a well-posed and superlinearly convergent sequence of iterates in the basic SQP method even for simple problems of nonlinear programming.}
\end{Remark}

\begin{Example}[\bf failure of solvability for subproblems in the basic SQP method]\label{ex041}{\rm Consider the one-dimensional nonlinear program
\begin{equation}\label{nlp}
\min-\frac{1}{2}x^2+\frac{1}{6}x^3\quad\mbox{subject to}\quad x\ge 0.
\end{equation}
By definition \eqref{mapS} of the solution map $S$ for this problem we have
\begin{equation*}
S(v,w)=\big\{(x,\lm)\in\R^2\;\big|\;v=-x+\sm x^2 +\lm\quad\mbox{and}\quad x+w\in N_{\R_-}(\lm)\big\},
\end{equation*}
where $(v,w)\in\R\times\R$. Set $(\ox,\olm):=(0,0)$ and observe that $(\ox,\olm)\in S(0,0)$. It is not hard to check that $\Lambda(\ox)=\{\olm\}$ and that $\ox$ is {\em not} a local minimizer for problem \eqref{nlp}; indeed, $x=2$ is the unique minimizer for this problem. We show nevertheless that the solution map $S$ enjoys the isolated calmness property at $((0,0),(\ox,\olm))$. Let us verify to this end that the KKT mapping $G$ from \eqref{GKKT} in the case of \eqref{nlp} is strongly metrically subregular at $((\ox,\olm),(0,0))$ by using the criterion in \eqref{isoc}. Observe that the mapping $G$ for this framework reduces to
\begin{equation*}
G(x,\lm)=\left[\begin{array}{c}
-x+\sm x^2 +\lm\\
-x
\end{array}
\right]+\left[\begin{array}{c}
0\\
N_{\R_-}(\lm)
\end{array}
\right],
\end{equation*}
and its graphical derivative is easily calculated by
\begin{equation*}
DG\big((\ox,\olm),(0,0)\big)(\xi,\eta)=\left[\begin{array}{cc}
-1&1\\
-1&0
\end{array}
\right]\left[\begin{array}{c}
\xi\\
\eta
\end{array}
\right]+\left[\begin{array}{c}
0\\
DN_{\R_-}(\olm, \ox)(\eta)
\end{array}
\right].
\end{equation*}
Since $DN_{\R_-}(\olm,\ox)(\eta)=N_{\R_-}(\eta)$, we arrive at
\begin{equation*}
(0,0)\in DG\big((\ox,\olm),(0,0)\big)(\xi,\eta)\iff\xi=\eta\quad\mbox{and}\quad\xi\in N_{\R_-}(\eta)
\end{equation*}
from which it clearly follows that $\xi=\eta=0$. Then \eqref{isoc} tells us that the mapping $G$ is strongly metrically subregular at $\big((\ox,\olm),(0,0)\big)$.

Next we are going to show that the generalized equation \eqref{kktsub4} has no solution for all $(u,\mu)\in\B_{1/2}(\ox,\olm)$ with $u\ne\ox$.
To furnish it, pick $(u,\mu)\in\B_{1/2}(\ox,\olm)$ with $u\ne\ox$ and observe that for \eqref{nlp} the generalized equation in question is represented by
\begin{equation*}
\left[\begin{array}{c}
0\\
0
\end{array}
\right]\in
\left[\begin{array}{c}
-u+\sm u^2+\mu\\
-u
\end{array}
\right]+
\left[\begin{array}{cc}
-1+u&1\\
-1&0
\end{array}
\right]
\left[\begin{array}{c}
x-u\\
\lm-\mu
\end{array}
\right]
+\left[\begin{array}{c}
0\\
N_{\R_-}(\lm)
\end{array}
\right].
\end{equation*}
This leads us to the relationships
\begin{equation*}
\lm=\sm u^2-x(u-1),\;x\lm=0,\;\;x\ge 0,\;\;\lm\le 0.
\end{equation*}
From the second condition above we deduce that either $x=0$ or $\lm=0$. If the former holds, it follows from the first equation that
$\lm=\sm u^2>0$, a contradiction. If the latter holds, we get $x=\dfrac{u^2}{u-1}<0$, a contradiction as well. This justifies that the KKT system associated with the subproblems of the basic SQP method for \eqref{nlp}\ has no solution for such a pair $(u,\mu)$.

Since we have $\Lambda(\ox)=\{\olm\}$, the basic constraint qualification \eqref{rcq} is satisfied. This tells us that if subproblems \eqref{subs} corresponding to \eqref{nlp} have local optimal solutions associated with $(u,\mu)$, then we would end up with having a solution
for the generalized equation \eqref{kktsub4}, which is not possible. This demonstrates the failure of solvability for subproblems in the basic SQP.}
\end{Example}

Let us mention that assuming {\em metric regularity} vs.\ strong metric subregularity of the mapping $G=S^{-1}$ from \eqref{GKKT} ensures the solvability of subproblems in the basic SQP method; cf.\ \cite[Theorem~6D.2]{dr}. However, the next theorem shows that in our framework the metric regularity yields the strong metric subregularity of the mapping $G$, which is equivalent to the isolated calmness of solution map $S$ from \eqref{mapS}. Recall to this end that the isolated calmness of $S$ alone does not ensure the existence of the SQP iterates as shown in Example~\ref{ex041}.

\begin{Theorem}[\bf strong metric subregularity from metric regularity for KKT systems]\label{mr-sms} Let $(\ox,\olm)$ be a solution to the KKT system \eqref{VS} under the standing assumptions {\rm(H1)} and {\rm(H2)}. If the mapping $G$ from \eqref{GKKT} is metrically regular around $((\ox,\olm),(0,0))$, then it is strongly metrically subregular at this point.
\end{Theorem}
{\bf Proof}. According to the coderivative criterion for metric regularity \eqref{cod-cr}, the mapping $G$ enjoys this property around $((\ox,\olm),(0,0))$ if and only if
\begin{equation}\label{mrg}
(0,0)\in D^*G\big((\ox,\olm),(0,0)\big)(w,u)\implies w=0,\,u=0.
\end{equation}
It is not hard to see that for any $(w,u)\in\R^n\times\R^m$ the coderivative of $G$ is calculated by
\begin{equation*}
D^*G\big((\ox,\olm),(0,0)\big)(w,u)=\left[\begin{array}{c}
\nabla_{xx}^2L(\ox,\olm)w-\nabla f(\ox)^*u\\
\nabla f(\ox)w+D^*N_\Th^{-1}\big(\olm,f(\ox)\big)(u)\\
\end{array}
\right].
\end{equation*}
As follows from \cite[Theorem~4G.1]{dr}, $G$ is strongly metrically subregular at $((\ox,\olm),(0,0))$ if and only if we have the implication
\begin{equation}\label{smsg}
(0,0)\in DG\big((\ox,\olm),(0,0)\big)(w,u)\implies w=0,\,u=0.
\end{equation}
Furthermore, it is easy to calculate as in \eqref{isoc} that
\begin{equation*}
DG\big((\ox,\olm),(0,0)\big)(w,u)=\left[\begin{array}{c}
\nabla_{xx}^2 L(\ox,\olm)w+\nabla f(\ox)^*u\\
-\nabla f(\ox)w+DN_\Th^{-1}\big(\olm,\Phi(\ox)\big)(u)\\
\end{array}
\right]
\end{equation*}
for all $(w,u)\in\R^n\times\R^m$. We get from Theorem~\ref{proto}(ii) under our standing assumptions {\rm(H1)} and {\rm(H2)} that the normal cone mapping $N_\Th$ is proto-differentiable at $f(\ox)$ for $\olm$. The latter together with \cite[Theorem~13.57]{rw} gives us the {\em derivative-coderivative inclusion}
\begin{equation}\label{der-coder}
DN_\Th\big(f(\ox),\olm\big)(\xi)\subset D^*N_\Th\big(f(\ox),\olm)(\xi)\quad\mbox{for all}\;\;\xi\in\R^m.
\end{equation}
Pick now $(w,u)\in\R^n\times\R^m$ from the left-hand side of implication \eqref{smsg} and deduce from the derivative-coderivative inclusion \eqref{der-coder} that $(w,-u)$ satisfies the relation on the left-hand side of implication \eqref{mrg}. Hence $(w,u)=(0,0)$, which yields \eqref{smsg}, and thus ensures the strong metric subregularity of $G$ at $((\ox,\olm),(0,0))$. $\h$\vspace*{0.05in}

A similar result to the one in Theorem~\ref{mr-sms} was obtained in \cite[Lemma~4F.8]{dr} for NLPs by using the graphical derivative characterization of the isolated calmness in \eqref{isoc}. It was extended in \cite[Corollary~25]{dsz} to ${\cal C}^2$-cone reducible constrained optimization problems by using a delicate topological result by Fusek \cite{fu}. The obtained Theorem~\ref{mr-sms} further extends the aforementioned results to the general case of parabolically regular constrained problems by employing a completely different device based on the coderivative criterion, the derivative-coderivative inclusion, and our recent second-order developments on parabolic regularity.

\section{Conclusions}\label{sect5}

This paper provides first applications of the novel theory of parabolic regularity in second-order variational analysis to characterizations of robust stability properties of KKT systems in general nonpolyhedral problems of constrained optimization with the subsequent usage of them to the development and justification of numerical algorithms. In particular, we derive complete second-order characterizations of robust isolated calmness for KKT systems with respect to arbitrary (not just canonical) perturbations. The obtained results are applied to developing the basic SQP method of solving constrained optimization problems governed by parabolically regular sets. This approach allows us to fully justify well-posedness of the basic SQP method (including solvability and stability of the corresponding subproblems) and its primal-dual superlinear convergence for parabolically regular constrained problems. Most of the obtained results are new not only in the general framework under consideration, but also for any conventional setting of nonpolyhedral constrained optimization.

Our future research aims at further developments of the presented approach to make it powerful in applications to more involved and efficient versions of the generic SQP algorithm as well as in other primal-dual methods of constrained optimization.

\end{document}